\magnification = \magstep1

\centerline{\bf Lectures on Maximal Monotone Operators} \medskip \centerline
{\bf  R. R. Phelps} 
\medskip \noindent Dept. Math. GN--50, Univ. of Wash., Seattle WA 98195;
phelps@math.washington.edu

\smallskip \noindent  (Lectures given at Prague/Paseky Summer School, Czech
Republic, August 15--28, 1993.)

\bigskip \noindent {\bf Introduction.}  These lectures will focus on those
properties of maximal monotone operators which are valid in arbitrary real Banach
spaces.  Most applications (to nonlinear partial differential equations,
optimization, calculus of variations, etc.) take place in reflexive
spaces, in part because several key properties have only been shown to
hold in such spaces.  (See, for instance, [De], [Pa--Sb] and [Ze].)  We will
generally isolate the reflexivity hypothesis, hoping that by doing so, it will
eventually be possible to decide their validity without that hypothesis.  In
Section 1 we define maximal monotone operators and prove some of their main
elementary properties.  Section 2 is devoted to the prototypical class of
subdifferentials of convex functions.  Gossez's subclass of monotone
operators of type (D) is examined in Section 3 and Section 4 gives a
brief treatment of another subclass, the locally maximal monotone
operators.  We have attempted to keep the exposition self--contained, using
only standard tools of elementary functional analysis.  One exception is the
application of the Brouwer fixed--point theorem in the proof of
the Debrunner--Flor theorem (Lemma 1.7).

\bigskip \noindent {\bf 1.  Monotone Operators.}

\proclaim Definition 1.1. A set-valued map  $T$  from a Banach space $E$  into
the subsets of its dual  $E^*$  is said to be a {\it monotone operator} provided
$$\langle x^* - y^*, x - y\rangle \ge 0 \quad \forall x, y\in E\   {\rm and}\ 
x^*\in T(x),\  y^*\in T(y).$$  We do not require that  $T(x)$  be nonempty.  The
{\it domain} (or {\it effective domain}) of  $T$  is the set $D(T) = \{x \in
E\colon \ T(x) \neq \emptyset\}$.

\medskip\noindent{\bf  Examples 1.2.}

(a)  The simplest examples of such operators are linear and single--valued.  For
instance, if $H$  is a real Hilbert space and  $T\colon H \to H^* \equiv H$  is a
linear map,  then  $T$  is monotone if and only if it is, in the usual sense, a
{\it positive operator}:  $\langle T(x), x\rangle\ge 0$  for all  $x$. 

  \medskip

(b)  Let $D$ be a nonempty subset of the real numbers ${\bf R}$.  A function 
$\varphi \colon D \to {\bf  R}^* \equiv {\bf R}$ defines a monotone operator if
and only if $\varphi$ is monotone nondecreasing in the usual sense:  That is,  
$$[\varphi(t_2) - \varphi(t_1)]\cdot(t_2 - t_1) \ge 0\  \forall\  t_1, t_2\in
D\quad {\rm iff}\quad \varphi(t_1) \le \varphi(t_2)\  {\rm whenever}\  t_1 <
t_2.$$ 

\medskip (c)  Examples of  {\it set--valued} monotone functions from ${\bf  R}$
to subsets of ${\bf  R}$ are easy to exhibit:  For instance, let $\varphi (x) =
0$ if $x < 0$, $\varphi (x) =
1$ if $x > 0$ and let $\varphi (0)$ be {\it any} subset of $[0, 1]$.

 \medskip (d)  Here is an important single--valued but nonlinear example:  Let 
$f$  be a continuous real--valued function on $E$ which is Gateaux differentiable
(that is, for each $x \in E$ the limit $$\displaystyle df(x)(y) = \lim_{t\to 0}
{f(x + ty) - f(x)\over t},\qquad y \in E$$ exists and is a bounded linear
functional of $y$).  Such a function {\it $f$ is convex if and only if the
mapping $x \to df(x)$ is monotone}.  Indeed, suppose that $f$ is convex; then for
$0 < t < 1$, convexity implies that $$\displaystyle {f(x + t(y - x)) - f(x) \over
t} \leq {(1 - t)f(x) + tf(y) - f(x) \over t} = f(y) - f(x).$$ It follows that
$df(x)(y - x) \leq f(y) - f(x)$, for any $x, y \in E$.  Thus, if $x, y \in E$ and
$x^* = df(x)$, $y^* = df(y)$, then $$\langle x^*, y - x\rangle \le f(y) -
f(x)\hbox{ and }-\langle y^*, y - x\rangle = \langle y^*, x - y\rangle \le f(x) -
f(y);$$ now add these two inequalities.  A proof of the converse may be found in
[Ph, p.17]. \medskip

(e)  The next example arises in fixed-point theory.  Let  $C$  be a bounded
closed convex nonempty subset of Hilbert space  $H$  and let  $U$  be a
(generally nonlinear) nonexpansive map of  $C$  into itself:  $\Vert U(x) -
U(y)\Vert\le\Vert x - y\Vert$  for all  $x$, $y\in C$.  Let  $I$  denote the
identity map in  $H$;  then  $T = I - U$  is monotone, with  $D(T) = C$.  Indeed,
for all  $x$, $y\in C$,   $$\langle T(x) - T(y), x - y\rangle  = \langle x - y
-(U(x) - U(y)), x - y\rangle =$$  $$= \Vert x - y\Vert^2 - \langle U(x) - U(y), x
- y\rangle \ge \Vert x - y\Vert^2 - \Vert U(x) - U(y)\Vert\cdot\Vert x - y\Vert
\ge 0.$$

\noindent Note that $0$ is in the range of $T$ if and only if $U$ has a fixed
point in  $C$;  this hints at the importance for applications of studying the
{\it ranges} of monotone operators. \medskip

(f)  Again, in Hilbert space, let $C$ be a nonempty closed convex set and let $P$
be the metric projection of $H$ onto $C$; that is, $P(x)$ is the unique element of
$C$ which satisfies $\|x - P(x)\| = \inf \{ \|x - y\|\colon y \in C\}$.  We first
prove the fundamental fact that the mapping $P$ satisfies (in fact, it is
characterized by) the following {\it variational inequality}:  For all $x \in
H$,  $$\langle x - P(x), z - P(x) \rangle \leq 0 \quad {\rm for\ all}\ z \in
C.\eqno (1.1)$$  Indeed, if $z \in C$ and  $0 < t < 1$, then $z_t \equiv tz + (1
- t)P(x) \in C$ and hence $\|x - P(x)\| \leq \|x - z_t\| = \|(x - P(x)) - t(z -
P(x))\|$.  Squaring both sides of this inequality, expanding and then cancelling
$\|x - P(x)\|^2$ on both sides yields $$0 \leq -2t\langle x - P(x), z -
P(x)\rangle + t^2\|z - P(x)\|^2.$$  If we then divide by $t$ and take the limit
as $t \to 0$ we obtain (1.1).  Moreover, if $y \in H$ and we write down (1.1)
again, using $y$ in place of $x$, then take $z = P(y)$ in the first equation, $z
= P(x)$ in the second one and add the two, we obtain
 $$\langle x - y, P(x) - P(y) \rangle \geq \|P(x) - P(y)\|^2 \quad {\rm for\
all}\ x, y \in H,\eqno (1.2)$$  which shows that $P$ is in monotone in a very
strong sense.  Note that $P$ is an example of a nonexpansive mapping in the sense
of the previous example:  We always have the inequality $\langle x - y, P(x) -
P(y)\rangle \leq \|x - y\|\cdot \|P(x) - P(y)\|$, so combined with the
monotonicity inequality above, we have $\|P(x) - P(y)\| \leq \|x - y\|$ for all
$x, y \in H$.

\medskip (g)  Here is a fundamental example of a set--valued monotone mapping,
the {\it duality mapping} from $E$ into $2^{E^*}$.  For any $x \in E$ define
$$J(x) = \{x^*\in E^*\colon \langle x^*, x\rangle  = \Vert x^*\Vert\cdot \Vert
x\Vert\hbox{ and } \Vert x^*\Vert = \Vert x\Vert\}.$$ By the Hahn--Banach
theorem, $J(x)$ is nonempty for each $x$, so $D(J) = E$.  Suppose that $x^* \in
J(x)$ and $y^* \in J(y)$.  Then $$\langle x^* - y^*, x - y\rangle = \|x^*\|^2 -
\langle x^*, y\rangle - \langle y^*, x\rangle + \|y^*\|^2 \geq \|x^*\|^2 -
\|x^*\|\cdot \|y\| - \|y^*\|\cdot \|x\| + \|y^*\|^2$$ $$ = \|x^*\|^2 -
2\|x^*\|\cdot \|y^*\| + \|y^*\|^2 = (\|x^*\| - \|y^*\|)^2,$$ so this is also
monotone in a rather strong sense.

\bigskip In order to define {\it maximal} monotone operators we must consider
their graphs.

\proclaim Definition 1.3.  A subset  $G$  of  $E \times E^*$  is said to be {\it
monotone} provided  $\langle x^* - y^*, x - y\rangle  \ge 0 $ whenever  $(x,
x^*)$, $(y, y^*)\in G$.  A set--valued mapping $T\colon E \to  2^{E^*}$  is a
monotone operator if and only if its graph   $$G(T) = \{(x, x^*)\in E \times
E^*\colon  x^*\in T(x)\}$$ is a monotone set. A monotone set is said to be {\it
maximal monotone} if it is maximal in the family of monotone subsets of  $E \times
E^*$,  ordered by inclusion. An element $(x, x^*) \in E \times E^*$ is said to be
{\it monotonically related} to the subset $G$ provided $$ \langle x^* - y^*, x -
y\rangle \geq 0 \quad {\rm for\ all}\ (y, y^*) \in G.$$ We say that a monotone
{\it operator}  $T$  is maximal monotone provided its graph is a maximal monotone
set.  

\medskip \noindent The most frequently used form of the definition of
maximality of $T$ is the following condition:    Whenever $(x, x^*) \in E \times
E^*$ is monotonically related to $G(T)$, then $x \in D(T)$ and $x^*\in T(x)$.   

\medskip\noindent There is an obvious one-to-one correspondence between monotone
sets and monotone operators.  An easy application of Zorn's lemma shows that {\it
every monotone operator}  $T$  {\it can be extended to a maximal monotone
operator}  $\overline  T$ ,  {\it in the sense that}  $G(T) \subset G(\overline 
T)$.

\proclaim Definition 1.4.  If $T\colon E \to  2^{E^*}$  is a monotone operator,
its {\it inverse} $T^{-1}$ is the set valued mapping from $E^*$ to $2^E$ defined
by $T^{-1}(x^*) = \{x \in E\colon\ x^* \in T(x)\}$.  Obviously,  $$G(T^{-1}) =
\{(x^*, x) \in E^* \times E\colon\ x^* \in T(x)\},$$ which is (within a
permutation) the same as the monotone set $G(T)$.  In particular,then, $T^{-1}$ is
maximal monotone if and only if $T$ is maximal monotone.

\vfill \eject

\noindent {\bf Examples 1.5.}  

(a) The monotone mapping $\varphi $ defined in Example 1.2 (c) is maximal if and
only if $\varphi (0) = [0, 1]$.  More generally, it is easily seen that a
monotone nondecreasing function $\varphi$ on ${\bf R}$ is maximal monotone if and
only if $\varphi (x) = [\varphi (x^-), \varphi (x^+)]$ for each $x \in {\bf R}$
(where, for instance, $\varphi (x^-) \equiv \lim_{t \to x^-} \varphi (t)$).

\medskip (b)  Any positive linear operator $T$ on Hilbert space is maximal
monotone.  Indeed, suppose $(x, x^*) \in H \times H$ is monotonically related to
$G(T)$.  Then, for any $z \in H$ and $\lambda > 0$ we have $$0 \leq \langle T(x
\pm \lambda z) - x^*, (x \pm \lambda z) - x\rangle = \pm \lambda \langle T(x)
\pm \lambda T(z) - x^*, z\rangle$$ $$ = \pm \lambda \langle T(x) - x^*, z\rangle +
\lambda^2 \langle T(z), z\rangle .$$  Dividing by $\lambda$ and then letting
$\lambda \to 0$ shows that $\langle
 T(x) - x^*, z \rangle = 0$ for all $z \in H$, hence that $x^* = T(x)$.

\medskip \noindent
{\bf Exercise 1.6.}  Prove that if  $T$  is maximal monotone,
then  $T(x)$  is a convex set, for every  $x\in E$.

\bigskip A major goal of these lectures is to examine the consequences of maximal
monotonicity with regards to questions of convexity of $D(T)$ and $R(T)$ (or
convexity of their interiors or closures).  Other basic questions involve
conditions under which $D(T) = E$ or $R(T) = E^*$, or whether these sets might be
dense in $E$ or $E^*$, respectively.  

\bigskip  Perhaps the most fundamental result concerning monotone operators is
the extension theorem of Debrunner--Flor [D--F].  An easy consequence of the
following version of their theorem states that {\it if $T$ is maximal
monotone and if its range is contained in a weak* compact convex set $C$, then
$D(T) = E$.}  (Given $x_0 \in E$, let $\phi$ (below) be the constant mapping
$\phi (x^*) = x_0$ for all $x^* \in C$, by the lemma, there exists $x^*_0 \in C$
such that $\{(x_0, x_0^*)\}  \cup G(T)$ is monotone; by maximality, $x_0 \in
D(T)$.)

\bigskip \proclaim  Lemma 1.7 (Debrunner--Flor).  Suppose  that  $C $ is a weak*
compact  convex subset of  $E^*$, that $\phi \colon\ C \to E$ is weak* to norm
continuous and that  $M \subset E \times C$  is a  monotone set.  Then there
exists  $x_0^* \in C $ such that  $\{(\phi (x_0^*), x_0^*)\}  \cup M$  is a
monotone  set.

\noindent {\bf Proof.}   For each element  $(y, y^*) \in M$  let $$U(y, y^*) =
\{x^* \in C\colon \langle x^* - y^*, \phi (x^*) - y\rangle < 0\}.$$ 
Since $x^* \to \langle x^* - y^*, \phi (x^*) - y\rangle $ is weak* continuous on
the bounded set $C$, each of these sets is 
relatively weak* open.  If the conclusion of the lemma fails,
then $C = \bigcup \{ U(y, y^*)\colon\ (y, y^*) \in M\}$.  By compactness,
there must exist $(y_1, y_1^*), \ (y_2, y_2^*), \ ..., \ (y_n, y_n^*)$ in $M$
such that $C = \bigcup_{i = 1}^n \{ U(y_i, y_i^*)\}$.  Let $\beta_1, \beta_2,
\dots , \beta_n$ be a partition of unity subordinate to this cover of $C$; that
is, each $\beta_i$ is weak* continuous on $C$, $0 \leq \beta_i \leq 1$, $\Sigma
\beta_i = 1$ and $\{x^* \in C\colon\ \beta_i(x^*) > 0\} \subset U(y_i, y^*_i)$
for each $i$.  Let $K = {\rm co}\{y_i^*\} \subset C$ and define the weak*
continuous map $p$ of $K$ into itself by
$$p(x^*) = \Sigma\beta_i(x^*) y_i^*, \quad x^* \in K.$$
Note that $K$ is a finite dimensional compact convex set which (since the weak*
topology is the same as the norm topology in finite dimensional spaces) is
homeomorphic to a finite dimensional ball.  Thus, the Brouwer fixed--point
theorem is applicable.  (See [Fr] for several proofs of the latter.)  It follows
that there exists $z^* \in K$ such that $p(z^*) = z^*$. We therefore have 
$$0 =
\langle p(z^*) - z^*, \Sigma \beta_j(z^*)(y_j - \phi (z^*))\rangle = \langle 
\Sigma \beta_i(z^*)(y_i^* - z^*), \Sigma \beta_j(z^*)(y_j -\phi (z^*))\rangle
=$$ $$= \Sigma _{i, j} \beta_i(z^*) \beta_j(z^*)\langle y_i^* - z^*, y_j - \phi
(z^*)\rangle.$$ 
Define $\alpha_{ij} = \langle y_i^* - z^*, y_j - \phi (z^*)\rangle$. 
It is straightforward to verify that $$\alpha_{ij} + \alpha_{ji} = \alpha_{ii} +
\alpha_{jj} + \langle y_i^* - y_j^*, y_j - y_i\rangle \leq \alpha_{ii} +
\alpha_{jj},$$ the inequality following from the monotonicity of $M$.  Next, to
simplify notation, let $\beta_i(z^*) = \beta_i$.  Note that for all $i, j$,
$$\beta_i \beta_j \alpha_{ij} + \beta_j \beta_i \alpha_{ji} = \beta_i \beta_j
({\alpha_{ij} + \alpha_{ji}  \over 2}) + \beta_j \beta_i ({\alpha_{ij} +
\alpha_{ji}  \over 2}) .$$
It follows that 
$$0 = \Sigma \beta_i \beta_j \alpha_{ij}= \Sigma \beta_i \beta_j ({\alpha_{ij}
+ \alpha_{ji}  \over 2}) \leq \Sigma \beta_i \beta_j ({\alpha_{ii} +
\alpha_{jj}  \over 2}).$$ 
We claim that this inequality implies that $\beta_i \beta_j = 0$ for all
$i, j$.  Indeed, for every pair $i, j$ such that
$\beta_i \beta_j > 0$ we must have $z^* \in U(y_i, y^*_i) \cap U(y_j, y_j^*)$,
hence both $\alpha_{ii} < 0$ and $\alpha_{jj} < 0$ so that $\Sigma \beta_i
\beta_j ({\alpha_{ii} + \alpha_{jj}  \over 2}) < 0$, a contradiction.  We
conclude that $\beta_i \equiv \beta_i(z^*) = 0$ for all $i$, an impossibility,
since $\Sigma \beta_i(z^*) = 1 $.

 \proclaim Definition 1.8.  A set--valued mapping $T\colon E \to 2^{E^*}$ is said
to be {\it locally bounded} at the point $x \in E$ provided there exists a
neighborhood $U$ of $x$ such that $T(U)$ is a bounded set.

\medskip \noindent Note that this does not require that the point $x$ actually be
in $D(T)$.  Thus, it is true (but not interesting) that $T$ is locally bounded at
each point of $E\setminus \overline {D(T)}$. 

\medskip \noindent There are at least four proofs (all based on the
Baire category theorem) that a maximal monotone operator $T$ is locally bounded
at the interior points of $D(T)$.  (See the discussion in [Ph, Sec. 2].)  The
first one was by Rockafellar [Ro$_2$], who showed that much of it can be carried
out in locally convex spaces.  His proof is the longest one, but it has the
great advantage of simultaneously proving that the interior of $D(T)$ is
convex.  By specializing his proof to Banach spaces (below), it is possible to
shorten it considerably.  We first make some simple general observations about an
arbitrary subset $D \subset E$, its convex hull ${\rm co}D$ and its interior
${\rm int}\thinspace D$.   
\smallskip
\noindent (i) Always, $D \subset {\rm co} D$, so ${\rm int}\thinspace D \subset
{\rm int}({\rm co} D)$.  Suppose it is shown that ${\rm int}({\rm co} D) \subset
D$.  Since ${\rm int}({\rm co} D)$ is open, it is therefore a subset of ${\rm
int}\ D$ and hence ${\rm int}({\rm co} D) = {\rm int}\thinspace D$, showing that
${\rm int}\ D$ is convex.  (None of these assertions assume that ${\rm
int}\thinspace D$ is nonempty!)   

 \smallskip \noindent (ii)  However, {\it if ${\rm int}({\rm co}D)$ is nonempty
and contained in $D$ (so that ${\rm int}({\rm co}D) = {\rm int}\thinspace D$),
then $\overline D = \overline {{\rm int}({\rm co}D)}$}, hence is convex.  [This
follows from the fact that for any convex set $C$ with interior we have $C \subset
\overline {{\rm int}\thinspace C}$, hence $D \subset {\rm co}D \subset
\overline{{\rm int}({\rm co}D)}$ therefore $\overline D \subset \overline {{\rm
int}({\rm co}D)} = \overline {{\rm int}\thinspace D}\subset \overline D$.]

\smallskip \noindent (iii)  Another useful elementary fact is the following: 
{\it If $\{C_n\}$ is an increasing sequence of closed convex sets
having nonempty interior, then ${\rm int}\thinspace \bigcup C_n \subset \bigcup
{\rm int}\thinspace C_n$ }(and, in fact, equality holds).  [Here is a proof, the
reader might have a simpler one:  If $x \in {\rm int}\thinspace \bigcup C_n $,
then it is in some $C_n$, hence in the closure of  $\bigcup {\rm int}\thinspace
C_n$.  The latter is an open convex set, so if it doesn't contain $x$, there must
exist a closed half--space $H$ supporting its closure at $x$.  Thus, ${\rm
int}\thinspace C_n \subset H$ for each $n$, therefore $C_n = \overline {{\rm
int}\thinspace
 C_n} \subset H$.  This leads to the contradiction that $x \in {\rm int}\thinspace
\bigcup C_n \subset {\rm int}\thinspace H$.]

\bigskip \noindent {\bf Notation}.  We will denote the closed unit ball in $E$
[resp. in $E^*$] by $B$ [resp. $B^*$].  Thus, if $r > 0$, say, then
$$rB^* = \{x^* \in E^*\colon \|x^* \| \leq r\}.$$

\proclaim Theorem 1.9 (Rockafellar).  Suppose that $T$ is maximal monotone and
that ${\rm int\ co} D(T)$ is nonempty.  Then ${\rm int}\thinspace D(T) = {\rm
int\ co} D(T)$ (so ${\rm int}\thinspace D(T)$ is convex) and $T$ is locally
bounded at each point of ${\rm int}\thinspace D(T)$.  Moreover, $\overline{ D(T)}
= \overline {{\rm int}\thinspace D(T)}$, hence it is also convex.

\noindent {\bf Proof.}   Let $C = {\rm int\ co} D(T)$.  For each $n \geq 1$ let
$$S_n = \{x \in nB\colon\ T(x) \cap nB^* \neq \emptyset\}.$$ Then $S_n \subset
S_{n+1}$ and $D(T) = \bigcup S_n \subset \bigcup {\rm co} S_n$.  Since the
$S_n$'s are increasing, this last union is convex, so it contains ${\rm co} D(T)$
and therefore contains $C$.  As an open subset of a Banach space, C has the Baire
property and therefore there exists an integer $n_0$ such that the closure
(relative to $C$) of $C \cap {\rm co}S_{n}$ has nonempty interior for all $n \geq
n_0$.  In particular, the larger set ${\rm int} (\overline {{\rm co}}S_{n})$
is nonempty for each such $n$.  We have 
$${\rm
int\ co} D(T) \subset {\rm int} \bigcup_{n \geq n_0} \overline {{\rm
co}}S_n \subset \bigcup_{n \geq n_0}{\rm int}(\overline {{\rm co}}S_n),$$ the last
inclusion being a special case of the elementary fact (iii) described above.  We
will show two things; first, that $T$ is locally bounded at any point in each
set ${\rm int} (\overline {{\rm co}}S_n)$ ($n \geq n_0$) and second, that each
such point is in $D(T)$ (which, as noted in (i) above, will imply that ${\rm
int}\ D(T)$ is convex).  For the first step, then, suppose that $x_0 \in {\rm int}
(\overline {{\rm co}}S_n)$ (for a fixed $n \geq n_0$).  Assume without loss of
generality that $n$ is sufficiently large that $R(T) \cap nB^* \neq \emptyset$
and for each $m \geq n$ let  
$$M_m = \{(u, u^*) \in E \times E^*\colon\  u \in D(T)\ {\rm and}\ u^* \in T(u)
\cap mB^* \};$$
this is a nonempty monotone subset of $E \times mB^*$.  For each $m \geq n$ and
$x \in E$ define 
  $$A_m(x) = \{x^* \in E^*\colon\ \langle x^* - u^*, x
- u \rangle \geq 0\quad  \forall u \in D(T)\ {\rm and}\ u^* \in T(u) \cap mB^*
\}.$$  
Since $A_m(x)$ is the set of all $x^*$ such that $(x, x^*)$ is monotonically
related to $M_m$, the Debrunner--Flor Lemma 1.7 guarantees that it is nonempty. 
For each such $m$ and every $x \in E$ we have $A_{m+1}(x) \subset A_m(x)$ and
$T(x) \subset A_m(x)$.  Moreover, as the intersection of weak*--closed
half--spaces, each $A_m(x)$ is weak* closed.  Suppose that $u^* \in nB^*$; then
$S_n \subset \{x \in E\colon\ |\langle u^*, x\rangle | \leq n^2\}$ and therefore
$\overline {{\rm co}}S_n$ is contained in the same set.  By hypothesis, there
exists $\epsilon > 0$ such that $x_0 + 2\epsilon B \subset \overline {{\rm
co}}S_n$.  Now, choose any $x \in x_0 + \epsilon B$ and $x^* \in A_n(x)$.  For
all $u \in S_n$ and $u^* \in T(u) \cap nB^*$ we must have $$\langle x^*, u - x
\rangle \leq \langle u^*, u - x\rangle \leq 2n^2,$$ which shows that $S_n \subset
\{u \in E\colon\ \langle x^*, u - x\rangle \leq 2n^2\}$ hence $\overline {{\rm
co}} S_n$ is contained in this same (closed and convex) set.  Thus, $$x_0 +
2\epsilon B \subset \overline { {\rm co}} S_n \subset \{u \in E\colon\ \langle
x^*, u - x\rangle \leq 2n^2\}.$$   Suppose that $\|v \| \leq \epsilon$, so that
$x + v \in x_0 + 2\epsilon B \subset \overline { {\rm co}} S_n$.  We then have
$\langle x^*, v\rangle = \langle x^*, (x + v) - x\rangle \leq 2n^2$.  Thus,
$\epsilon \|x^*\| = \sup \{\langle x^*, v\rangle \colon\ \|v \| \leq \epsilon\}
\leq 2n^2$, which implies that $\|x^*\| \leq 2n^2/\epsilon$.  We have shown,
then, that if $x \in x_0 + \epsilon B$ and $x^* \in A_n(x)$, then $x^* \in
{2n^2\over \epsilon}B^*$.  Since $$T(x_0 + \epsilon B) = \bigcup \{T(x)\colon\ x
\in x_0 + \epsilon B\} \subset \bigcup \{A_n(x)\colon \ x \in x_0 + \epsilon B\}
\subset (2n^2/ \epsilon)B^*,$$ we see that $T$ is locally bounded at $x_0$.
(Note that this part of the proof did not require that $T$ be maximal.)

\smallskip \noindent To show that $x_0 \in D(T)$, note that $A_n(x_0) \subset
(2n^2/ \epsilon)B^*$ and is weak* closed, hence is weak* compact.  For $m
\geq n$ we have  $A_m(x_0) \subset A_n(x_0)$ and therefore the sequence
$\{A_m(x_0)\}$ is a decreasing family of nonempty weak* compact sets.  Let $x_0^*
\in \bigcap_{m \geq n} A_m(x_0)$.  If $u \in D(T)$ and $u^* \in T(u)$, then
$\|u^*\| \leq m$ for some $m \geq n$.  By definition of $A_m(x_0)$ we have
$\langle x_0^* - u^*, x_0 - u\rangle \geq 0$.  By maximality of $T$, this implies
that $x_0^* \in T(x_0)$ and therefore $x_0 \in D(T)$.

\proclaim Corollary 1.10.  Suppose that $E$ is reflexive and that $T\colon E \to
2^{E^*}$ is maximal monotone.  Then ${\rm int}\thinspace  R(T)$ is
convex.  If ${\rm int}\thinspace  R(T)$ is nonempty, then $\overline
{R(T)}$ is convex.

\noindent {\bf Proof.}  Since $T^{-1}$ is a maximal monotone operator from $E^*$
to $2^E$, by the previous theorem, ${\rm int}\thinspace  D(T^{-1})
\equiv {\rm int}\thinspace  R(T)$ is convex.

\bigskip \noindent  That the first part of this corollary can fail in a
nonreflexive space will be shown by Example 2.21 (below).

\proclaim Definition 1.11.  A subset $A \subset E$ (not necessarily convex) which
contains the origin is said to be {\it absorbing} if  $E = \cup\{\lambda A: 
\lambda > 0\}$.  Equivalently,  $A$  is absorbing if for each  $x\in E$  there
exists  $t > 0$  such that  $tx \in A$.  A point  $x\in A$  is called an {\it
absorbing point of} $A$  if the translate  $A - x$  is absorbing.

\medskip \noindent It is obvious that any interior point of a set is an absorbing
point.  If  $A_1$  is the union of the unit sphere and  $\{0\}$,  then  $A_1$  is
absorbing, even though it has empty interior.  A proof of the following theorem
of J. Borwein and S. Fitzpatrick [B--F] may be found in [Ph].

\proclaim Theorem 1.12.  Suppose that $T\colon E \to  2^{E^*}$  is monotone and
that  $x\in D(T)$.  If  $x$  is an absorbing point of  $D(T)$ (in particular, if
$x\in \hbox{\rm int }D(T)$), then $T$ is locally bounded at $x$.

\noindent Note that the foregoing result does not require that  $D(T)$  be convex
nor that $T$ be maximal.  There are trivial examples which show that  $0$  can be
an absorbing point of  $D(T)$  but not an interior point (for instance, let  $T$ 
be the restriction of the duality mapping $J$ to the set  $A_1$  defined above).
Even if  $D(T)$  {\it is} convex and  $T$  is maximal monotone,  $D(T)$ can have
empty interior, as shown by the following example.  (In this example,  $T$  is an
unbounded linear operator, hence it is not locally bounded at any point and
therefore  $D(T)$  has no absorbing points.) \medskip

\proclaim Example 1.13.  In the Hilbert space  $\ell^2$  let  $D = \{x =
(x_n)\in \ell^2\colon\  (2^nx_n)\in \ell^2\}$  and define  $Tx = (2^nx_n)$,  $x\in
D$.  Then  $D(T) = D$  is a proper dense linear subspace of  $\ell^2$  and  $T$ 
is a positive operator, hence -- by Example 1.5 (b) -- it is maximal monotone.

\medskip    {\it It is conceivable that for a maximal monotone  $T$,  any
absorbing point of  $D(T)$  is actually an interior point.}  That this is true if
$\overline {D(T)}$ is assumed to be convex is shown by combining Theorem 1.12
with the following result.

\proclaim Theorem 1.14 (Libor Vesel\'y). Suppose that $T$ is maximal monotone and
that $\overline {D(T)}$ is convex.  If $x \in \overline {D(T)}$ and $T$ is
locally bounded at $x$, then $x \in {\rm int}\ D(T) $.

\noindent {\bf Proof.}  The first step doesn't use the convexity hypothesis:  {\it
Suppose that $T$ is a maximal monotone operator which is locally bounded at the
point $x \in \overline {D(T)}$; then $x \in D(T)$.}  Indeed, by hypothesis, there
exists a neighborhood $U $ of $x$ such that $T(U)$ is a bounded set.  Choose a
sequence $\{x_n\} \subset D(T)\cap U$ such that $x_n \to x$ and choose $x_n^* \in
T(x_n)$.  By weak* compactness of bounded subsets of $E^*$ there exists a subnet
of $\{(x_n, x_n^*)\}$ -- call it $\{(x_\alpha, x_\alpha^*)\}$ -- and $x^* \in
 E^*$ such that $x_\alpha^* \to x^*$ (weak*).  It follows that for all $(y, y^*)
\in G(T)$,
 $$\langle x^* - y^*, x - y\rangle = \lim_\alpha \ \langle x_\alpha^* - y^*,
x_\alpha - y\rangle \geq 0;$$ by maximal monotonicity, $x^* \in T(x)$, so $x \in
D(T)$.  Next, if $x$ is in the boundary of the closed convex set $\overline
{D(T)}$, then $T$ is {\it not} locally bounded at $x$:  Suppose there were a
neighborhood $U$ of $x$ such that $T(U)$ were bounded.  By the Bishop--Phelps
theorem there would exist a point $z \in U \cap \overline {D(T)}$ and a nonzero
element $w^* \in E^*$ which supported $\overline {D(T)}$ at $z$; that is,
$\langle w^*, z\rangle = \sup \ \langle z^*, D(T) \rangle$.  Now, $T$ would also
be locally bounded at $z \in U$, so by the first step, $z \in D(T)$ and we could
choose $z^* \in T(z)$.  For any $(y, y^*) \in G(T)$ and any $\lambda \geq 0$ we
would have $$\langle z^* + \lambda w^* - y^*, z - y\rangle = \langle z^* - y^*, z
- y\rangle + \lambda \langle w^*, z - y\rangle \geq 0.$$  By the maximality of $T$
this would imply that  $z^* + \lambda w^* \in T(z)$ for each $\lambda \geq 0$,
which shows that $T(z)$ is not bounded, a contradiction.  Since $x \notin {\rm
bdry}\ \overline {D(T)}$, it must be in ${\rm int}\thinspace  \overline {D(T)}$.
By local boundedness, we can choose an open set $U$ such that $x \in U \subset
{\rm int}\thinspace  \overline {D(T)}$ and $T(U)$ is bounded.  Thus, $T$ is
locally bounded at {\it every} point of $U$, which -- by the first step proved
above -- implies that $U \subset D(T)$ and therefore $x \in {\rm int}\thinspace 
D(T)$. \medskip \noindent Note that this result, combined with Theorem 1.9,
implies that for any maximal monotone $T$, {\it if ${\rm int}\thinspace D(T)$ is
nonempty, then it is precisely the set of points in $\overline {D(T)}$ where
$T$ is locally bounded.}  It remains open as to what happens if ${\rm
int}\thinspace D(T)$ is empty and $\overline {D(T)}$ is not convex.

\proclaim Exercise 1.15.  Prove that if  $T$  is maximal monotone, then
for all $x\in \hbox{ int }D(T)$ the set $T(x)$ is weak* compact and convex.

\proclaim Definition 1.16.  Let $X$ and $Y$ be Hausdorff spaces and
suppose that $T\colon X \to 2^{Y}$ is a set--valued mapping.  We say that $T$ is
{\it upper semicontinuous at the point} $x \in X$ if the following holds:  For
every open set $U \subset Y$ such that $T(x) \subset U$ there exists an open
subset $V$ of $X$ such that $x \in V$ and $T(V) \subset U$.  Upper
semicontinuity on a {\it set} is defined in the obvious way.

\proclaim Exercise 1.17.  Prove that if $T\colon E \to 2^{E^*}$ is maximal
monotone, then it is norm--to--weak* upper semicontinuous on ${\rm int}\thinspace
D(T)$.

\bigskip The following ``fixed--point'' result, which will be useful to us in
Section 3, illustrates both the utility of the upper semicontinuity property
as well as the strength of the Debrunner--Flor Lemma 1.7.  It is a special
case of Theorem 4 of [Bro].

\proclaim Lemma 1.18.  Suppose that $E$ is a reflexive Banach space and that
$K$ is a nonempty compact convex subset of $E$.  Let $R\colon\ K \to 2^K$ be
an upper semicontinuous mapping such that $R(u)$ is nonempty, closed and
convex, for each $u \in K$.  Then there exists $u_0 \in K$ such that $u_0 \in
R(u_0)$.

\noindent {\bf Proof.}  Suppose there were no such point; then 
$0 \notin u - R(u)$ for each $u \in
K$.  By the separation theorem applied to the compact convex set $u - R(u)$
and $0$, for each $u$ there would exist $x^* \in E^*$, $\|x^*\| = 1$, and
$\delta > 0$ such that $\langle x^*, v\rangle > \delta$ for each $v \in u -
R(u)$.  For each $x^* \in E^*$, define $$W(x^*) = \{u \in K\colon\ \langle x^*,
v\rangle > 0 \quad \forall\ v \in u - R(u)\}.$$   For each $\|x^*\| = 1$ let 
$U(x^*) = \{v \in E\colon\ \langle x^*, v\rangle >0\}$, so $u \in W(x^*)$ if
and only if $u \in K$ and $u - R(u) \subset U(x^*)$.   Now, if $u \in K$, then
our supposition implies that there exists $\|x^*\| =
1$ and $\delta > 0$ such that $u - R(u) + \delta B \subset U(x^*)$.  Upper
semicontinuity of $R$ (hence of $-R$) at $u$ implies that for some $0 < \epsilon <
\delta /2$ we will have $-R(y) \subset -R(u) + {\delta \over 2}B$ whenever $y \in
(u + \epsilon B) \cap K$.  It follows that for all such $y$, we have $$y - R(y)
\subset u + {\delta \over 2}B -R(u) + {\delta \over 2}B \subset U(x^*),$$
that is, $(u + \epsilon B) \cap K \subset W(x^*)$.  This shows that every point
$u \in K$ is in the interior of some $W(x^*)$, so that the sets $\{{\rm
int}\thinspace W(x^*)\}$ form an open cover of $K$.  As in the proof of Lemma
1.7, there is  a finite subcover $\{{\rm int}\thinspace W(x_j^*)\}_{j = 1}^n$ of
$K$ and a continuous partition of unity $\{ \beta_1, \beta_2, \dots, \beta_n\}$
subordinate to this cover.  Define $$r(x) = \sum \beta_j(x) x_j^*, \quad x \in
K.$$ This is a continuous map from $K$ into $E^*$ and for all $u \in K$ and $v
\in u - R(u)$ we have $\langle r(u), v\rangle = \sum \beta_j(u) \langle x_j^*,
v\rangle > 0$ (since $\beta_j(u) > 0$ implies that $u \in W(x_j^*)$ hence that
$\langle x_j^*, v\rangle > 0$).  In Lemma 1.7 let $\phi = -r$, reverse the roles
of $E$ and $E^*$ (reflexivity permits this) and let $M$ be the monotone set $ K
\times \{0\}$ to obtain $u_0 \in K$ such that $-\langle r(u_0), u_0 - u\rangle
\geq 0$ for all $u \in K$.  In particular, this is true if $u = v_0$ where $v_0$
is any element of $ R(u_0)$; that is, the element $v \equiv u_0 - v_0 \in u_0 -
R(u_0)$ satisfies $\langle r(u_0), v\rangle \leq 0$, in contradiction to the fact
that $\langle r(u_0), v\rangle > 0$, completing the proof.

\bigskip \medskip

\noindent {\bf 2.  Subdifferentials of Convex Functions}

\medskip This section introduces what is perhaps the most basic class of maximal
monotone operators.  We assume that the reader has some familiarity with
real--valued convex functions.  (A good reference for their elementary
properties is [R--V].)  Their use in optimization and convex analysis is most
simply handled by introducing the seeming complication of admitting {\it extended}
real-valued functions, that is, functions with values in ${\bf R}
\cup\{\infty\}$.

\proclaim Definition 2.1.  Let  $X$  be a Hausdorff space and
let  $f\colon X \to  {\bf R} \cup  \{\infty \}$.  The {\it effective domain}
of  $f$  is the set  $\hbox{dom}(f) = \{x\in X\colon f(x) < \infty \}$. 
Recall that  $f$  is {\it lower semicontinuous} provided  $\{x\in X\colon f(x) \le 
r\}$  is closed in $X$  for every  $r\in {\bf R}$.  This is equivalent to saying
that the {\it epigraph} of  $f$ $$\hbox{epi}(f) = \{(x, r)\in X\times {\bf
R}\colon  r \ge f(x)\} $$
is closed in  $X \times{\bf R}$.  Equivalently,  $f$  is lower
semicontinuous provided
$$f(x) \le  \lim \inf f(x_\alpha)$$  
whenever  $x\in X$  and  $(x_\alpha)$  is a net in  $X$  converging
to  $x$. We say that  $f$  is {\it proper} if  $\hbox{dom}(f) \neq 
\emptyset$.  

\medskip \noindent Note that if $f$ is defined on a Banach space $E$ and is
convex, then so is  $\hbox{dom}(f)$. Also,  {\it a function}  $f$  {\it is convex
if and only if} $\hbox{epi}(f)$  {\it is convex}.  This last fact is important; it
implies that certain properties of lower
semicontinuous convex functions can be deduced
from properties of these (rather special) closed
convex subsets of  $E\times {\bf R}$.  One can view this as saying
that the study of lower semicontinuous convex
functions is a special case of the study of closed
convex sets.

\proclaim Examples 2.2.
(a)  Let  $C$  be a nonempty convex
subset of  $E$; then the {\it indicator function}  $\delta_C$, 
defined by  $\delta_C(x) = 0$  if  $x\in C, \quad  = \infty$   otherwise, is
a proper convex function which is lower
semicontinuous if and only if  $C$  is closed.

\noindent
This example is one reason for introducing
extended real-valued functions, since it makes
it possible to deduce certain properties of a
closed convex set from properties of its lower
semicontinuous convex indicator function.  Thus,
one can cite this example to support the view that the study
of closed convex sets is a special case of the
study of lower semicontinuous convex functions. 
It's all a matter of which approach is more
convenient, the geometrical or the analytical. 
It is useful to be able to switch easily from one to
the other.

\medskip (b)  Let  $A$  be any nonempty subset of  $E^*$   such
that the weak* closed convex hull of  $A$  is not all
of  $E^*$  (or, more simply, let  $A$  be a weak* closed
convex proper subset of $E^*$) and define the {\it support function}
$\sigma_A$ of $A$ by   $$\sigma_A(x) = \hbox{sup}\{\langle x^*, x\rangle : x^*\in
A\},\quad   x\in E.$$   This is easily seen to be a proper lower semicontinuous
convex function. 

\medskip (c)  If  $f$  is a continuous convex function defined
on a nonempty closed convex set  $C$, extend  $f$  to be 
$\infty$   at the points of  $E\setminus C$;  the resulting function is
a proper lower semicontinuous convex function.

\bigskip
The next proposition uses completeness of  $E$  to
describe a set where a lower semicontinuous convex function
is necessarily continuous.

\proclaim Proposition 2.3. Suppose that  $f$  is a proper
lower semicontinuous convex function on a Banach
space  $E$   and that  $D = \hbox{\rm int } \hbox{\rm dom}(f)$  is
nonempty;  then  $f$  is continuous on  $D$.

\noindent {\bf Proof.}  We need only show that  $f$  is locally
bounded in  $D$, since this implies that it is locally
Lipschitzian in  $D$ (see [Ph, Prop. 1.6]).  First, note that if  $f$  is
bounded above (by  $M$,  say) in   $B(x; \delta) \subset  D$  for some 
$\delta  > 0$,   then  it is bounded  below  in  $B(x; \delta)$. 
Indeed, if  $y$  is in  $B(x; \delta)$,  then so is  $2x - y$  and 
$$f(x) \le  {1\over2}[f(y) + f(2x - y)] \le 
{1\over2}[f(y) + M]$$
so  $f(y) \ge 2f(x) - M$  for all  $y\in B(x;\delta)$.  Thus, to
show that  $f$  is locally bounded in  $D$,  it suffices
to show that it is locally bounded {\it above} in  $D$.  For
each  $n \ge 1$,  let  $D_n = \{x\in D: f(x) \le 
n\}$. The sets  $D_n$
are closed and  $D = \cup D_n$; since  $D$  is a Baire space,
for some  $n$  we must have  $U \equiv  \hbox{ int }D_n$ nonempty. We
know that  $f$  is bounded above by  $n$  in  $U$;  without
loss of generality, we can assume that  $B(0;\delta) \subset  U$ 
for some  $\delta  > 0$.  If $y$ is in  $D$,  with  $y \neq  0$, then
there exists  $\mu  > 1$   such that   $z = \mu y\in D$   and hence
(letting   $0 < \lambda  = \mu^{-1} < 1)$,  the set  
$$V = \lambda z +(1 - \lambda)B(0;\delta) = y + (1 - \lambda )B(0;\delta)$$  
is a neighborhood of  $y$  in  $D$.  For any point  $v = (1 - \lambda )x +
\lambda z\in V \quad (\hbox{\rm where}\quad x \in B(0;\delta ))$  we have  
$$f(v) \le  (1 - \lambda )n + \lambda f(z),$$          
so  $f$  is bounded above in  $V$  and the proof is
complete.
\medskip

\noindent{\bf Examples 2.4.}

 (a)  The function  $f$  defined by  
$f(x)
= 1/x$  on  $(0, \infty )$,  $f(x) = \infty$   on $(-\infty , 0]$  shows that 
$f$  can be continuous at a boundary point  $x$  of 
$\hbox{dom}(f)$  where  $f(x) = \infty$.  (Recall that the
neighborhoods of  $\infty$   in  $(-\infty , \infty ]$  are all the sets $(a,
\infty ]$,  $a\in {\bf  R}$.)

\medskip (b)  Suppose that  $C$  is nonempty closed and
convex; then the lower semicontinuous convex indicator function $\delta_C$  is
continuous at  $x\in C$  if and only if  $x\in \hbox{int }C$.  Thus, {\it if 
$\hbox{\rm int }C = \emptyset $,  then  $\delta_C$  is not continuous at any point
of}  $C = \hbox{ dom}(\delta_C)$.

\proclaim Definition 2.5. Recall that if  $E$  is a Banach
space, then so is  $E\times {\bf R}$,  under any norm which
restricts to give the original topology on the subspace 
$E$,  for instance,  $\Vert (x, r)\Vert = \Vert x\Vert + \vert r\vert $.  Recall,
also, that  $(E\times{\bf R})^*$  can be identified with   $E^* \times{\bf
R}$,  using the pairing 
$$\langle (x^*, r^*), (x, r)\rangle  = \langle x^*,
x\rangle  + r^*\cdot r.$$ 

\noindent {\bf Remark 2.6.}  If a proper lower semicontinuous
convex function  $f$  is continuous at some point  $x_0
\in \hbox{dom}(f)$,  then   $\hbox{dom}(f)$  has nonempty interior and 
$\hbox{epi}(f)$  has nonempty interior in  $E\times {\bf R}$.  (Indeed,  
$f(x) = \infty $  outside of  $\hbox{dom}(f)$, so  $x_0$  cannot be a
boundary point of the latter.  Moreover, there
exists an open neighborhood  $U$  of  $x_0$  in  $\hbox{dom}(f)$  in
which  $f(x) < f(x_0) + 1$,  so the open product set  $U \times
\{r: r > f(x_0) + 1\}$  is contained in  $\hbox{epi}(f)$.)

\proclaim Definition 2.7.  If  $x\in
\hbox{dom}(f)$, define the {\it subdifferential mapping} $\partial f$ by  
$$\partial f(x) = \{x^*\in E^* : \langle x^*, y - x\rangle 
\le  f(y) - f(x) \hbox{ for all } y\in E\}$$
$$= \{x^*\in E^* : \langle x^*, y\rangle  \le  f(x + y) - f(x) \hbox{
for all } y\in E\},$$
while  $\partial f(x) = \emptyset$   if  $x\in
E\setminus\hbox{dom}(f)$.  {\it It may also be empty at points of }
$\hbox{dom}(f)$, as shown in the first example below.

\medskip \noindent  It is easy to see that $\partial f$ is a monotone operator: 
If $x^*\in \partial f(x)$ and $y^*\in \partial f(y)$, then  $$\langle x^*, y -
x\rangle \le f(y) - f(x)\hbox{ and }-\langle y^*, y - x\rangle = \langle y^*, x -
y\rangle \le f(x) - f(y);$$ now add these two inequalities.

\medskip \noindent  It is not obvious that $\partial f$ is {\it maximal} monotone;
in fact, it is not even obvious that it is nontrivial (i.e., that $D(\partial f)
\neq \emptyset$.)  Much of the rest of this section is devoted to
establishing these properties.

\proclaim Definition 2.8.  If  $x\in \hbox{dom}(f)$  we define its right--hand
directional derivative  $d^+f(x)$  by $$d^+f(x)(y) = \lim_{t\to 0^+} t^{-1}[f(x +
ty) - f(x)],\qquad     y \in E.$$

\noindent It follows from the convexity of $f$ that this limit always exists (see
[Ph]).  Note that  $d^+f(x)(y) = \infty $  if  $x + ty\in
E\setminus\hbox{dom}(f)$  for all  $t > 0$.  (It is also possible to have 
$d^+f(x)(y) = -\infty$;  consider, for instance,  $d^+f(0)(1)$  when  $f(x) =
-x^{1/2}$  for  $x \ge 0$, $  = \infty$  elsewhere.)  We have the following
important relationship:  For any point  $x\in \hbox{dom}(f)$,
$$x^*\in \partial f(x) \hbox{ if and only if } \langle x^*,
y\rangle  \le  d^+f(x)(y) \hbox{ for all } y\in E. $$
It follows from this that for the example given above ($f(x) = -\sqrt x$ for $x
\geq 0$), it must be true that $\partial f(0) = \emptyset$.  In the first
example below, one sees that it is possible to have $\partial f(x) = \emptyset$
for a dense set of points $x \in {\rm dom}(f)$.

\bigskip
\noindent{\bf Examples 2.9.}

(a) Let $C$ be the closed (in fact, compact)
convex subset of $\ell^2$ defined by 
$$C=\{x\in\ell^2:\vert
x_n\vert \le 2^{-n}, n=1, 2, 3,\ldots\}$$
 and define $f$ on $C$ by $f(x)=\Sigma [-(2^{-n}+x_n)^{1/2}].$
Since each of the functions $x \to -(2^{-n} + x_n)^{1/2}$ is continuous, convex
and bounded in absolute value by $2^{(-n+1)/2}$,  the
series converges uniformly, so  $f$  is continuous
and convex. We claim that $\partial f(x) = \emptyset$ for any $x\in C$ 
such that  $x_n > -2^{-n}$  for infinitely many  $n$. Indeed, let  ${\bf e}_n$ 
denote the  $n$-th unit vector in  $\ell^2$.  If  $x^*\in \partial f(x)$ (so
that, as noted above,  $x^* \le  d^+f(x)$),  then for all  $n$  such that $x_n >
-2^{-n}$, we have   $$-\Vert x^*\Vert  \le  \langle x^*, {\bf e}_n\rangle  \le 
d^+f(x)({\bf e}_n) = -(1/2)(2^{-n} + x_n)^{-1/2},$$
an impossibility which implies that  $\partial f(x) = \emptyset $.  Note
that if we make the usual extension (setting  $f(x) =
\infty $  for  $x\in \ell^2\setminus C)$,  then  $f$  is lower
semicontinuous, but not continuous at any point of  $C \ \ (= \hbox{ bdry
}C)$.
\medskip (b) Let $C$ be a nonempty closed convex subset of $E$; then for any
$x\in C$, the subdifferential $\partial\delta_C(x)$ of the indicator function
$\delta_C$ is the cone with vertex $0$ of all $x^*\in~E^*$ which ``support''
$C$ at $x$, that is, which satisfy $$\langle x^*, x\rangle=\hbox{ sup}\{\langle
x^*, y\rangle:y\in C\}\equiv\sigma_C(x^*).$$
(Indeed, $ x^*\in \partial \delta_C(x)$  if and only if  $\langle x^*,
y - x\rangle  \le  \delta_C(y) - \delta_C(x)$, while  $x^*$  attains its supremum
on  $C$  at $ x$  if and only if the left hand side of this latter
inequality is at most  $0$,  while the right side is
always greater or equal to  $0$.)

\medskip The following notion of an {\it approximate} subdifferential is useful
in many parts of convex analysis.

\proclaim Definition 2.10. Let  $f$  be a proper
convex lower semicontinuous function and suppose  $x\in \hbox{dom}(f)$. 
For any  $\epsilon  > 0$  define the $\epsilon $--{\it subdifferential} 
$\partial_\epsilon f(x)$  by       
$$\partial_\epsilon f(x) = \{x^*\colon \langle x^*, y\rangle  \le  f(x +
y) - f(x) + \epsilon \hbox{  for all } y\in E\}.$$

\noindent It follows easily from the definition that $\partial_\epsilon
f(x)$ convex and weak* closed.  That fact that it is nonempty for every $x \in
{\rm dom}(f)$ follows from the convexity of ${\rm epi}(f)$ and the separation
theorem (in $E \times {\bf R}$) (see [Ph; Prop. 3.14]).

\medskip  The basic maximality technique which was used in proving the
Bishop--Phelps theorem [Ph] was applied in $E \times {\bf R}$  by Br\o ndsted
and Rockafellar to prove the following fundamental lemma, which shows, among
other things, that $\partial f$ is nontrivial.

\proclaim Lemma 2.11.  Suppose that  $f$  is a convex proper lower semicontinuous
function on the Banach space  $E$.  Then given any  point $x_0 \in \hbox{\rm
dom}(f)$,  $\epsilon > 0$, $\lambda > 0$ and any functional $ x_0^* \in
\partial_\epsilon f(x_0)$, there exist  $x \in \hbox{\rm dom}(f)$ 
and  $x^*\in  E^*$   such that   $$x^* \in \partial f(x), \quad    \Vert x -
x_0\Vert\le \epsilon/\lambda \quad  \hbox{  and  }\quad   \Vert x^* -
x_0^*\Vert\le  \lambda.$$   In particular, the domain of  $\partial f$  is dense
in  $\hbox{\rm dom}(f)$, so $\overline {D(\partial f)} = \overline
{{\rm dom}(f)}$ is convex.

\medskip  The next result can be looked at in two ways.  On the one hand, it is
simply a verification that (under a certain hypothesis), the subdifferential
operation is additive.  (We use the usual vector sum of sets.)  On the other
hand, it will be seen later as a verification in a special case that the sum of
maximal monotone operators is again maximal monotone.

\proclaim Theorem 2.12. Suppose that  $f$  and  $g$  are convex
proper lower semicontinuous functions on the
Banach space  $E$  and that there is a point in 
$\hbox{\rm dom}(f) \cap\hbox{\rm dom}(g)$  where one of them, say  $f$,  is
continuous.  Then
$$\partial (f + g)(x) = \partial f(x) + \partial
g(x),  \qquad         x\in D(\partial f) \cap D(\partial g).$$

{\bf Remark.}  It is immediate from the definitions
that for  $x\in\hbox{ dom}(f + g)$  (which is identical to 
$\hbox{dom}(f) \cap \hbox{dom}(g))$,  one must have  
$$\partial f(x) + \partial g(x) \subset  \partial (f + g)(x).$$

\noindent This inclusion can be proper.  To see this, let  $E
= {\bf R}^2$,  let  $f$ denote the indicator function  $\delta_C$  and let $g =
\delta_L$,  where  $C$  is the epigraph of the quadratic function  $y = x^2$
and  $L$  is the $x$-axis.  Obviously,  $C$  and  $L$  intersect only at the
origin $0$  and it is easily verified that  $\partial f(0) = {\bf R^-e}$, 
where  {\bf e}  is the vector  $(0, 1)$,  and  $\partial g(0) = {\bf Re}$, 
while  $$\partial (f + g)(0) = {\bf R}^2 \neq  \partial f(0) +
\partial g(0).$$
\noindent {\bf Proof.}  Suppose that  $x_0^* \in \partial (f + g)(x_0)$.  In order
to simplify the argument, we can replace  $f$  and  $g$ 
by the functions
$$f_1(x) = f(x + x_0) - f(x_0) -\langle x_0^* , x\rangle   \hbox{ and }  g_1(x) =
g(x + x_0) - g(x_0), \quad   x\in E;$$
it is readily verified from the definitions that if 
 $x_0^* ~\in ~\partial (f + g)(x_0)$,    then  $0\in \partial
(f_1 + g_1)(0)$  and if  $0\in \partial f_1(0) + \partial g_1(0)$,
then  $x_0^* \in \partial f(x_0) + \partial g(x_0)$. Without loss of
generality, then, we assume that  $x_0 = 0$,  $x_0^*  = 0$,  $f(0) = 0$ and 
$g(0) = 0$.  We want to conclude that  $0$  is in the sum  $\partial f(0) +
\partial g(0)$,  under the hypothesis that   $0\in \partial (f +
g)(0)$.  This last means that
$$(f + g)(x) \ge (f + g)(0) = 0 \hbox{  for all  }
x\in E.\eqno(2.1)$$
We now apply the separation theorem in  $E\times {\bf R}$  to the
two closed convex sets  $C_1 = \hbox{epi}(f)$  and $C_2 = \{(x,
r)\colon r \le  -g(x)\}$;  this is possible because  $f$  has a
point of continuity in  $\hbox{dom}(f) \cap \hbox{dom}(g)$  and hence --
recall Remark 2.6 -- $C_1$  has nonempty interior. 
 Moreover, it follows from (2.1) that  $C_2$  misses
the interior of  $C_1 = \{(x, r)\colon r > f(x)\}$.  Since  $(0,
0)$  is common to both sets, it is contained in any
separating hyperplane.  Thus, there exists a
functional  $(x^*, r^*)\in E^*  \times {\bf R}$,  $(0, 0) \neq  (x^*,
r^*)$,  such that
$$\langle x^*, x\rangle  + r^*\cdot r \ge 0 \hbox{ if } r \ge f(x)  \hbox{ and
}  \langle x^*, x\rangle  + r^*\cdot r \le  0 \hbox{ if } r \le  -g(x).$$
Since  $1 > f(0) = 0$  we see immediately that $ r^* \ge 0$. 
To see that  $r^* \neq  0$,  (that is, that the separating
hyperplane is not ``vertical''), we argue by
contradiction:  If  $r^* = 0$,  then we must have $ x^* \neq 
0$; also  $\langle x^*, x\rangle  \ge 0$  for all  $x\in
\hbox{ dom}(f)$  and  $\langle x^*, x\rangle  \le  0$  for all  $x\in
\hbox{ dom}(g)$. This says that  $x^*$  separates these two sets.  This is
impossible; by the continuity hypothesis, their intersection
contains an interior point of  $\hbox{dom}(f)$.  Without
loss of generality, then, we can assume that  $r^* =
1$  and hence, for any $x\in E$,
$$\langle -x^*, x - 0\rangle  \le  f(x) - f(0)  \hbox{ and } \langle x^*, x -
0\rangle  \le  g(x) - g(0),$$
that is,  $0 = -x^* + x^*\in \partial f(0) + \partial g(0)$,  which
completes the proof.

\medskip \noindent The next two lemmas lead easily to S. Simons' recent
proof [Si$_{1, 3}$] of Rockafellar's maximal monotonicity theorem for
subdifferentials.

\proclaim Lemma 2.13.  Suppose that $f$ is a lower semicontinuous proper
convex function on $E$.  If $\alpha, \beta > 0,\ x_0\in E$ and $f(x_0) <  {\rm
inf}_E f +  \alpha \beta,$ then there exist $x \in E$ and $x^* \in \partial
f(x)$  such that $\|x - x_0\| < \beta $ and $\|x^*\| < \alpha.$

\noindent {\bf Proof.}  Choose 
$\epsilon > 0$ such that  $f(x_0) - \hbox {\rm inf}_E f <
\epsilon < \alpha \beta$ and then choose $\lambda$ such that $\epsilon /\beta <
\lambda < \alpha.$  It follows that $0 \in \partial_\epsilon f(x_0)$ so by
the Br\o ndsted--Rockafellar Lemma 2.11, there exist $x \in {\rm dom}(f)$ and $x^*
\in \partial f(x)$ such that $\|x^*\| \leq \lambda < \alpha$ and $\|x - x_0\| \leq
\epsilon/\lambda < \beta$.

\proclaim Lemma 2.14.  With $f$ as in the previous lemma, suppose that $x \in E$
(not necessarily in ${\rm dom}(f)$) and that 
${\rm inf}_E f <  f(x).$  Then there exist $z 
\in {\rm dom}(f)$ and $z^* \in \partial f(z)$ such that $$f(z) < f(x) \ {\rm
and}\ \langle z^*, x - z\rangle > 0.$$

\noindent {\bf Proof.}  Fix $\lambda \in {\bf R}$ such that
${\rm inf}_E f < \lambda < f(x)$ and let $$K = \hbox {\rm
sup}_{y \in E, y \not= x} {\lambda  - f(y)
\over \|y - x\|}.$$   
We first show that $0 < K < \infty$.  To that end, let $F = \{y \in
E\colon f(y) \leq \lambda\}$, so $F$ is closed, nonempty and $x \notin F.$ 
Since $\hbox {\rm dom}(f) \not= \emptyset$, one can apply  the separation
theorem in $E \times {\bf R}$ to find $u^* \in E^*$ and $r \in {\bf R}$ such
that  $f \geq u^* + r.$  Suppose that $y\in E$ and that $y \not= x$.  If $y
\in F$, then $$\lambda - f(y) \leq \lambda - \langle u^*, y\rangle - r \leq
|\lambda - \langle u^*, x\rangle - r| + \langle u^*,
 x - y\rangle$$ hence $${\lambda - f(y) \over \|y - x\|} \leq {|\lambda -
\langle u^*, x\rangle - r| \over \hbox {\rm dist}(x, F)} + \|u^*\|.$$ If $y
\notin F$, then ${\lambda - f(y) \over \|y - x\|} < 0$.  In either case, there
is an upper bound for ${\lambda - f(y) \over \|y - x\|} $, so $K < \infty$.  To
see that  $K > 0$, pick any $y \in E$ such that $f(y) < \lambda$.  Since
$\lambda  < f(x)$, we have  $ y \not= x$ and $K \geq {\lambda - f(y) \over
\|y - x\|} > 0$.

Suppose, now, that $0 < \epsilon < 1$, so that $(1 - \epsilon)K < K$ and
hence, by definition of $K$, there exists $x_0 \in E$ such that $x_0 \not=
x$ and $${\lambda - f(x_0) \over \|x_0 - x\|} > (1 - \epsilon)K.$$
For $z \in E$, let $N(z) = K\|z -x\|$; we have shown that
$(1 - \epsilon) N(x_0) + f(x_0) <
\lambda$, that is, $(N + f)(x_0) < \lambda + \epsilon
N(y)$.  We claim that $\lambda \leq {\rm inf}_E (N + f)$  Indeed, if $z = x$, then
we have $\lambda < f(x) = (N + f)(z)$, while if $z \not= x$, then
${\lambda - f(z) \over
 \|z - x\|} \leq K$,  from which it follows that $\lambda \leq (N + f)(z)$. 
Thus, we have shown that there is a point $x_0 \in E$,  $x_0 \not= x$, such that
 $$(N + f)(x_0) < \hbox{\rm inf}_E (N + f) +
\epsilon K \|x_0 - x\|.$$
We now apply Lemma 2.13 to $N + f$,  with $\beta = \|x_0 - x \|$ and $\alpha =
\epsilon K.$  Thus, there exists $z \in {\rm dom}(N + f) \equiv {\rm dom}(f)$ and
$w^* \in \partial (N + f)(z)$ such that $\|z - x_0\| < \|x - x_0\|$ and $\|w^*\|
< \epsilon K.$  It follows that $\|z - x\| > 0$.  From the sum formula (Theorem
2.12),  $$\partial (N + f)(z) = \partial N(z) + \partial f(z),$$ 
so there exist $y^* \in \partial
N(z)$ and $z^* \in \partial f(z)$ such that $w^* = y^* + z^*$.  Since $y^* \in
\partial N(z)$, we must have $\langle y^*, z - x\rangle \geq N(z) - N(x) =
K\|z - x\|$.  Thus $$\langle z^*, x - z\rangle = \langle y^*, z - x\rangle
+ \langle w^*, x - z\rangle \geq K\|z - x\| - \|w^*\|\cdot  \|x - z\| > $$
$$ >(1 - \epsilon)K\|z - x\| > 0.$$ Since $z^*
\in \partial f(z)$, we have $f(x) \geq f(z) + \langle z^*, x - z\rangle >
f(z)$, which completes the proof.

\proclaim Theorem 2.15 (Rockafellar).  If $f$ is a proper lower semicontinuous
convex function on a Banach space $E$, then its subdifferential $\partial f$ is a
maximal monotone operator.

\noindent {\bf Proof.}  Suppose that $x \in E$, that $x^* \in E^*$ and that $x^*
\notin \partial f(x)$.  Thus, $0 ~\notin ~\partial (f - x^*)(x)$, which implies
that ${\rm inf}_E(f - x^*) < (f - x^*)(x)$.  By Lemma 2.14 there exists $z \in {\rm
dom}(f - x^*) \equiv {\rm dom}(f)$ and $z^* \in \partial (f - x^*)(z)$ such that
$\langle z^*, z - x\rangle < 0$.  Thus, there exists $y^* \in \partial f(z)$
such that $z^* = y^* - x^*$, so that $\langle y^* - x^*, z - x\rangle < 0.$ 

\proclaim Example 2.16.  In a Banach space  $E$  define 
$j(x) = (1/2)\Vert x\Vert^2$; this is clearly continuous and convex. The monotone
duality mapping $J$ is actually the subdifferential of $j$, and hence is maximal
monotone. 

\noindent {\bf Proof.}  It is readily computed that  $d^+j(x)(y) =
\Vert x\Vert\cdot d^+\Vert x\Vert(y)$.  If  $x = 0$, then $d^+j(0)(y) = 0$ for
all  $y$,  hence is linear and therefore  $\partial j(0) = \{0\}$. 
Suppose, then, that  $x \neq 0$.  We know (from the remark following Definition
2.8) that $x^*\in \partial j(x)$ if and only if $x^* \le d^+j(x)$,  that is, if
and only if  $\Vert x\Vert^{-1}x^* \le d^+\Vert x\Vert$,  which is equivalent to 
$y^* \equiv \Vert x\Vert^{-1}x^*\in \partial \Vert x\Vert$, that is, if and only
if $\langle y^*, y - x\rangle  \le \Vert y\Vert - \Vert x\Vert$ for all  $y\in
E$.  If, in this last inequality, we take  $y = x + z, \Vert z\Vert\le 1$  and
apply the triangle inequality, we conclude that  $\Vert y^*\Vert\le 1$.  If we
take  $y = 0$,  we conclude that  $\Vert x\Vert\le \langle y^*, x\rangle  \le
\Vert y^*\Vert\cdot \Vert x\Vert$,  so   $\Vert y^*\Vert = 1$  and  $\langle y^*,
x\rangle  = \Vert x\Vert$,  which is equivalent to what we want to prove.  The
converse is easy:  If  $\Vert y^*\Vert = 1$  and  $\langle y^*, x\rangle  = \Vert
x\Vert$,  then for all  $y$  in $E$,  we necessarily have  $\langle y^*, y -
x\rangle  \le\Vert y\Vert - \Vert x\Vert$,  so  $y^*\in \partial \Vert x\Vert$.

\medskip  There are other interesting and useful properties of the duality
mapping $J$.  For instance, it is immediate from the original definition that
it satisfies $J(-x) = -J(x)$ and $J(\lambda x) = \lambda
J(x)$ for $\lambda > 0$.  Since it is the subdifferential of
the function ${1\over 2}\|x\|^2$, it is not surprising that it reflects
properties of the norm.  

\proclaim Proposition 2.17. \smallskip (a)  The norm in $E$ is
Gateaux differentiable (at nonzero points) if and only if $J$ is single valued. 
\smallskip (b)  The mapping $J$ is ``one--to--one'' (that is, $J(x) \cap J(y) =
\emptyset$ whenever $x \neq y$) if and only if the norm in $E$ is {\it strictly
convex}; that is, $\|x + y\| < 2$ whenever $\|x\| = 1$, $\|y\| = 1$ and $x \neq
y$. \smallskip (c)  Surjectivity of $J$ is equivalent to reflexivity of $E$.

\noindent {\bf Proof.}  Parts (a) and (b) are straightforward exercises.  Part
(c) is
an easy consequence of R. C. James' deep result that a Banach space $E$ is
reflexive if (and only if) each functional in $E^*$ attains its supremum on the
unit ball of $E$ at some point.  Indeed, if $J$ is surjective and $x^* \in
E^*$, then there exists $x \in E$ such that $x^* \in J(x)$, hence $\langle x^*,
{x\over \|x\|}\rangle = \|x^*\| \equiv \sup\{ \langle x^*, y\rangle \colon \|y\|
\leq 1 \}$, showing $x^*$ attains its supremum on the unit ball at $x\over
\|x\|$.

\bigskip \proclaim Definition 2.18.  If $S$ and $T$ are monotone operators and
$x \in D(S) \cap D(T)$, we define 
$$(S + T)(x) = S(x) + T(x) \equiv \{ x^* + y^*\colon\ x^* \in S(x),\ y^* \in
T(x)\},$$ 
while $(S + T)(x) = \emptyset$ otherwise.

\noindent  It is immediate that $S + T$ is also a monotone operator and that --
by definition -- $D(S + T) = D(S) \cap D(T)$.

\medskip The following theorem of Rockafellar [Ro$_1$] is
basic to many applications of maximal monotone operators.

\proclaim Theorem 2.19 (Rockafellar).  Suppose that $E$ is reflexive, that $S$
and $T$ are maximal monotone operators on $E$ and that $D(T) \cap {\rm
int}\thinspace D(S) \neq \emptyset$.  Then $S + T$ is maximal.

\noindent We refer to [Ro$_1$] for the proof of this
theorem (see, also, [B-C-P]), which relies partly on the fact that any reflexive
space can be renormed so that both it and its dual norm are Gateaux
differentiable (at nonzero points) [D--G--Z].  Since maximal monotonicity does
not depend on which norm defines the topology of $E$, this guarantees that the
duality mapping $J$ can be assumed to be single--valued, one--to--one and onto, a
useful step in the proof.

\medskip \noindent The situation is wide open in arbitrary Banach spaces.

 \proclaim Problem
2.20.  Suppose that $E$ is a nonreflexive Banach space and that $S$ and $T$ are
maximal monotone operators such that $D(T) \cap {\rm int}\thinspace
D(S) \neq \emptyset$; is $S + T$ necessarily maximal?  What about the special
case when $S$ is the subdifferential of the indicator function $\delta_C$ of a
closed convex set $C$ for which ${\rm int}\thinspace C \cap D(T) \neq \emptyset$?

\noindent {\bf Remark.} {\it The answer to the first question is affirmative when
$D(S) = E = D(T)$}; this was pointed out to us by Martin Heisler.  One can use the
local boundedness of $S$ and $T$ and weak* compactness to prove that the graph of
$G(S + T)$ is closed in the norm $\times$ weak* topology in $E \times E^*$; Lemma
2.2 of [F-P$_1$] then applies to show that $S + T$ is maximal. 

\medskip \noindent  Rockafellar [Ro$_1$, Theorem 3] has shown that the answer to
the second question is affirmative for certain single--valued monotone
operators $T$.  

\medskip \noindent The answer to the first question is also affirmative whenever
$S = \partial f$
 and $T = \partial g$, where $f$ and $g$ are proper lower semicontinuous convex
functions; this is a consequence of Theorem 2.12 (that the subdifferential of the
sum of two convex functions is  the sum of their subdifferentials).  Indeed, if
$\hbox {\rm int}D(T) \neq \emptyset$, then [since $\hbox {\rm int dom}(\partial
g) \subset \hbox {\rm int dom}(g)$ and $g$  is continuous on $\hbox {\rm int
dom}(g)$] we conclude that $g$ is continous at some point of $D(S) \cap D(T)
\subset \hbox {\rm dom}(f) \cap {\rm dom}(g)$, so
Theorem 2.12 implies  that
$S + T = \partial (f + g)$ and Theorem 2.15 implies that the latter is
maximal monotone.  The remark following the statement of Theorem 2.12 shows that,
even in a two dimensional Banach space, maximality of a sum can fail if 
$D(S)\cap {\rm int}D(T) $ is empty; in that example, the graph of $\partial f
+ \partial g$ is a proper subset of the graph of the maximal monotone operator
$\partial (f + g)$.

\bigskip  The foregoing discussion is an example of how subdifferentials
fulfill their role as prototypes when considering general questions about
maximal monotone operators.  Whenever a property is valid for
subdifferentials in arbitrary Banach spaces there is some hope that it also holds
for all maximal monotone operators.  On the other hand, of course, if it fails for
subdifferentials on nonreflexive spaces, the situation is obviously hopeless. 
The following example illustrates this with respect to the first assertion in
Corollary 1.10 (that ${\rm int}\thinspace R(T)$ is convex when $E$ is reflexive).

 \proclaim Example 2.21 (Simon Fitzpatrick).  There exists a continuous convex
function $f$ on the Banach space $c_0$ such that the interior of $R(\partial f)$
is not convex.

\medskip\noindent {\bf Proof.}  With the usual supremum norm on $c_0$, define
$g(x) = \|x\|$ and $h(x) = \|x - e_1\|$, where $e_1 = (1, 0, 0, \dots)$, and
let $f = g + h$.  Since $g$ and $h$ are continous and convex, we have
$\partial f(x) = \partial g(x) + \partial h(x)$  for each $x \in c_0$.  It is
straightforward to compute that $\partial g(x)$ is either $B^*$ if $x = 0$, or is
contained in the set $F$ of all finitely nonzero sequences in $\ell_1$
otherwise.  It is also easy to see that $\partial h(x) = \partial g(x - e_1)$ for
each $x$.  It follows (letting $e_1^*$ denote the corresponding element of
$\ell_1$) that $\partial f(0) = -e_1^* + B^*$ and $\partial f(e_1) = e_1^* +
B^*$, while $\partial f(x)$ is contained in $F$ if $x \neq 0, e_1$.  Since ${\rm
int}\ R(\partial f) \supset {\rm int}B^* \pm e_1^*$, if it were convex, it would
contain $0 = {{1}\over {2}}(e_1^* - e_1^*)$ and hence a neighborhood $U$ of $0$. 
But, for any $\lambda > 0$, the element $(0, {\lambda\over 2^2}, {\lambda\over
2^3}, {\lambda\over 2^4},\dots)$ is not in $F$, has distance from $\pm e_1^*$ equal
to $1 + {\lambda\over 2} > 1$ and (for sufficiently small $\lambda$) is in $U$.

\bigskip Another place one looks for prototypical properties is in Hilbert space
(see, for instance, [Bre]), where things are made much easier by the fact that the
duality mapping $J$ is replaced by the identity mapping.

\medskip 
We conclude this section with a look at an additional property which
characterizes subdifferentials within the class of maximal monotone operators. 
Details may be found, for instance, in [Ph].

\proclaim Definition 2.22.  A set-valued map  $T\colon E \to  2^{E^*}$  is said
to be $n${\it -cyclically monotone} provided   $$\sum_{1\le k\le n} \langle
x_k^*, x_k - x_{k-1}\rangle  \ge 0$$   whenever  $n \ge 2$ and  $x_0, x_1, x_2,
\ldots x_n\in E$,  $x_n = x_0$, and  $x_k^*\in T(x_k)$,  $k = 1, 2, 3, \ldots,
n$.   We say that  $T$  is  {\it cyclically monotone}  if it is  $n$-cyclically
monotone for every $n$. Clearly, a 2-cyclically monotone operator is monotone.  

\noindent{\bf   Examples 2.23.}

(a)  The linear map in  ${\bf \rm R}^2$  defined by  $T(x_1, x_2) = (x_2, -x_1)$ 
is positive, hence maximal monotone, but it is not 3-cyclically monotone:  Look at
the points  $(1, 1)$,  $(0, 1)$ and  $(1, 0)$.
\medskip
(b)  Let  $f$  be a proper lower semicontinuous convex function; then  $\partial
f$  is cyclically monotone.
\medskip \noindent The final theorem of this section shows that this is the {\it
only} such example [Ro$_3$].

 \proclaim Theorem 2.24 (Rockafellar).  If $T\colon E \to
2^{E^*}$ is maximal monotone and cyclically monotone, with $D(T) \neq \emptyset$,
then there exists a proper convex lower semicontinuous function $f$ on $E$ such
that $T = \partial f$.

\bigskip
\noindent {\bf 3.  Gossez's monotone operators of type (D).}
\bigskip  In 1971, J.--P. Gossez [Go$_1$] introduced the class of monotone
operators of ``dense type'' in order to extend to nonreflexive spaces some of the
basic known results about maximal monotone operators on reflexive spaces.  He
subsequently modified his definition [Go$_{3,5}$] to the one given
below.  One first identifies a Banach space $E$ with its canonically
embedded image $\widehat E$ in $E^{**}$.  Having done this, it is
natural to consider the graph $G(T)$ of a monotone operator $T$ to
be a subset of $E^{**} \times E^*$.

\proclaim Definition 3.1.  A monotone operator
$T\colon E \to 2^{E^*}$ is said to of {\it type (D)} provided it satisfies the
following property:  If $(x^{**}, x^*) \in E^{**} \times E^*$ is monotonically
related to $G(T)$, then there exists a net $(x_\alpha, x^*_\alpha) \in G(T)$ such
that $x_\alpha \to x^{**}$ in the $\sigma(x^{**}, x^*)$ topology, $(x_\alpha)$ is
bounded and $x_\alpha^* \to x^*$ in norm.

\medskip \noindent In using this definition it is convenient to extend $T$ to a
mapping $\overline T\colon E^{**} \to 2^{E^*}$ as follows:
\smallskip \noindent
We let $\overline T$ be the map whose graph
$G(\overline T) \subset E^{**} \times E^*$ consists of all
elements $(x^{**}, x^*) \in E^{**} \times E^*$ which are monotonically related
to $G(T)$.  

\smallskip \noindent {\it
The map $\overline T$ need not be monotone} [Go$_4$].  However, {\it if $T$ is
monotone of type (D), then $\overline T$ is maximal monotone}.  Indeed, suppose
that $(x^{**}, x^*)$ is monotonically related to $G(\overline T)$; it is then
obviously monotonically related to $G(T)$, hence is in $G(\overline T)$.  It
only remains to show that $\overline T$ is monotone.  Suppose that $(x^{**},
x^*), (y^{**}, y^*) \in G(\overline T)$.  By hypothesis, there exists a net
$(x_\alpha, x^*_\alpha) \in G(T)$ as above, hence 
$$0 \leq \langle y^{**} - \hat x_\alpha, y^* - x^*_\alpha\rangle = \langle
y^{**} - x^{**}, y^* - x^*\rangle + \langle x^{**} - \hat x_\alpha, y^* -
x^*\rangle + \langle y^{**} - \hat x_\alpha, x^* - x^*_\alpha \rangle.$$
Taking limits, we see that $\langle
y^{**} - x^{**}, y^* - x^*\rangle \geq 0$.

\medskip \noindent A word about terminology:  When $T$ is monotone of type (D),
then any element of $G(\overline T)$ is the limit of a certain net in $G(T)$, so
we can consider the latter as being dense in the former, hence Gossez's earlier
``dense type'' and the later use of ``type (D)''.

\medskip
\noindent {\bf Examples 3.2.} 

(a)  Define $T\colon {\bf R}\backslash \{0\} \to {\bf R}$ by $T(x) = 0$ for
all $x \neq 0$.  It is easily verified that $G(\overline T) = {\bf R} \times
\{0 \} \subset {\bf R}^2$, hence is maximal monotone, while $T$ is monotone of
type (D) (but not maximal).
\medskip
(b)  Suppose that $E$ is reflexive. If $T$ is maximal monotone, then $\overline
T = T$, hence $T$ is trivially maximal monotone of type (D), that
is, {\it in reflexive spaces, the maximal monotone operators
coincide with the maximal monotone operators of type (D).}

\medskip
(c)  If $f$ is a proper lower semicontinuous convex function on $E$,
then $\partial f$ is maximal monotone of type (D) [Go$_1$].  (This is not
obvious.  In essence, it uses the first step of Rockafellar's original proof
[Ro$_3$] of the maximal monotonicity of $\partial f$.)

\medskip \noindent The fundamental fact about maximal monotone operators of
type (D) is contained in Theorem 3.6 (below).  In order to formulate it, we need
to recall a few facts about convex functions and their Fenchel duals.

\proclaim Definition 3.3.  If $f$ is a proper convex lower semicontinuous
function on $E$, then the function $f^*$ defined on $E^*$ by
$$f^*(x^*) = \sup \{\langle x^*, x\rangle - f(x)\}, \quad x^* \in E^*$$
is proper, convex and  weak* lower semicontinuous on $E^*$.  One defines
$f^{**} \equiv (f^*)^*$ on $E^{**}$ analogously.

\proclaim Exercises 3.4.  

(a) If we let $j(x) = {1\over 2}\|x\|^2$, then $j^*(x^*) =
{1\over 2}\|x^*\|^2$ and $j^{**}(x^{**}) = {1\over 2}\|x^{**}\|^2$.

\smallskip (b) If $f$ is convex proper and lower semicontinuous, and $\epsilon
\geq 0$, then $\partial_\epsilon f$ can be characterized as follows:  $x^* \in
\partial_\epsilon f(x)$ if and only if $f(x) + f^*(x^*) \leq \langle x^*,
x\rangle + \epsilon$.  (Here, we take $\partial_0 f \equiv \partial f$.)

\bigskip The following theorem, which is a very
special case of a result of F. Browder [Bro, Theorem 10], provides a crucial
step in proving Theorem 3.6 (below).

\proclaim Theorem 3.5. Suppose that $F$ is a finite dimensional Banach
space with unit ball $B$ and duality map $J$.  Fix $r > 0$ and let $A$
denote the restriction of $-J$ to $rB$.  If $T\colon\ rB \to 2^{F^*}$ is
monotone, then there exists $(x, x^*) \in G(A)$ such that $G(T) \cup \{(x,
x^*)\}$ is monotone, that is, there exists $x \in rB$ and $x^* \in -J(x)
\subset rB^*$ such that $\langle x^* - y^*, x - y\rangle \geq 0$ for all $(y,
y^*) \in G(T)$.

\noindent {\bf Proof.}  For each $\epsilon > 0$ and finite subset $G \subset
G(T)$, let
$$H_{\epsilon, G} = \{(x, x^*) \in G(A)\colon\ \langle x^* - y^*, x - y\rangle
\geq -\epsilon\quad \forall\ (y, y^*) \in G\}$$
and
$$H = \{(x, x^*) \in G(A)\colon (x, x^*)\ {\rm is\ monotonically\ related\ to}\
G(T)\}.$$
We will have proved the theorem if we show that the set $H$ is nonempty.  It is
straightforward to verify that $H = \bigcap \{H_{\epsilon, G}\colon\ \epsilon >
0,\ \ G\ {\rm finite},\ G\subset G(T) \}$.  Moreover, from the compactness of
$rB \times rB^*$ and the definition of $J$, it follows that each $H_{\epsilon,
G}$ is compact.  The intersection of any finite number of these sets contains a
set of the same form, so to show that the family has the finite intersection
property it suffices to show that each of them is nonempty; it will then
follow by compactness that $H$ is nonempty.  Fix $\epsilon$ and $G$.  To show
that $H_{\epsilon, G}$ is nonempty, we define, for each $x \in rB$,
$$S(x) = \{x^* \in rB^*\colon\ \langle x^* - y^*, x - y\rangle > -\epsilon
\quad \forall\ (y, y^*) \in G\}$$
and for $x^* \in rB^*$,
$$S^{-1}(x^*) = \{x \in rB\colon\ \langle x^* - y^*, x - y \rangle > -\epsilon
\quad \forall\ (y, y^*) \in G\}.$$
Each set $S(x)$ is a finite intersection of relatively open subsets of $rB^*$,
hence is relatively open and each set $S^{-1}(x^*)$ is convex.  Moreover, each
of the latter sets is nonempty:  This is an application of Lemma 1.7, using
the monotonicity of $G$, reversing the roles played there by $E$ and
$E^*$ and letting $\phi\colon\ rB \to rB^*$ be the constant map with value $x^*$. 
It follows that the sets $\{S(x)\}$ form an open cover of the compact set $rB^*$,
hence there exist $\{x_1, x_2,\dots, x_n\} \subset rB$ such that $rB^* =
\bigcup_{j = 1}^n S(x_j)$.  As in the proof of Lemma 1.7, there exists a partition
of unity $\{\beta_1, \beta_2, \dots, \beta_n\}$ subordinate to this covering. 
Define the continuous mapping $p\colon\ rB^* \to rB$ by $p(x^*) = \sum_{j = 1}^n
\beta_j(x^*) x_j$.  We claim that $x^* \in S(p(x^*))$ for all $x^* \in rB^*$. 
Indeed, for every $j$ such that $\beta_j(x^*) > 0$ we have $x^* \in S(x_j)$, that
is, $x_j \in S^{-1}(x^*)$.  Since the latter is convex and since $p(x^*)$ is a
convex combination of $x_j$'s, it follows that $p(x^*) \in S^{-1}(x^*)$, which is
equivalent to $x^* \in S(p(x^*))$.  Next, define the set--valued mapping $R\colon\
rB^* \to 2^{rB^*}$ by $R(x^*) = -J(p(x^*))$.  Since $p$ is continuous and $J$ is
upper semicontinuous (Exercise 1.17 and Example 2.16), $R$ is upper
semicontinuous.  Let $K = rB^*$; by Lemma 1.18 there exists $x_0^* \in rB^*$ such
that $x_0^* \in R(x_0^*)$, that is, $x_0^* \in -J(p(x_0^*))$ and, of course,
$x_0^* \in S(p(x_0^*))$.  Letting $x_0 = p(x_0^*)$, this means that $x_0^* \in
(-J)(x_0) \cap S(x_0) \subset H_{\epsilon, G}$, which completes the proof.

\bigskip \noindent Note that since $\partial j^*\colon\ E^* \to 2^{E^{**}},$
its inverse $(\partial j^*)^{-1}$ is a mapping from $E^{**}$ to $2^{E^*}$.

\proclaim Theorem 3.6 (Gossez).  If $T$ is maximal monotone of
type (D), then for all $\lambda > 0$, $R(\overline T + \lambda (\partial
j^*)^{-1}) = E^*$.

\noindent {\bf Proof.}  Since $T$ is of type (D) if and only if the
same is true of $\lambda^{-1} T - x^*$ for each $\lambda > 0$ and $x^*
\in E^*$, we need only show that $0 \in
R(\overline T + (\partial
j^*)^{-1})$.  Let $\cal F$ denote the directed family of
all finite dimensional subspaces $F \subset E$ such that $D(T) \cap F \neq
\emptyset$, partially ordered by inclusion.  For each such $F$, let
$i_F\colon F \to E$ denote the natural injection, with adjoint
$i_F^*\colon E^* \to F^*$ (the restriction mapping to $F$).  Suppose,
now that $F \in \cal F$ and that $r > 0$.  Apply Theorem 3.5 to $F$, $F^*$,
$i_F^*Ti_F$ and $K_r = \{x \in F\colon \|x \| \leq r\}$, as follows:  Let $G$
be the graph in $K_r \times F^*$ of the restriction to $K_r$ of the monotone
operator $i_F^*Ti_F$.  (We assume that $r$ is sufficiently large that
$K_r \cap D(T) \neq \emptyset$.)  Let $A\colon K_r \to F^*$ be
$-i_F^*Ji_F$.  Thus, there exists an element of $G(A)$ -- call it
$(x_{_{F,r}}, -x_{_{F,r}}^*)$ -- which is montonely related to $G$; that
is, $\|x_{_{F,r}}\| \leq r$, $x_{_{F,r}}^* \in i_F^*J(x_{_{F,r}})$ 
(this uses $Ji_F = J$ in $F$) and $\langle -x_{_{F,r}}^* - y^*,
x_{_{F,r}} - y \rangle \geq 0 $ whenever $\|y \| \leq r$ and $y^* \in
i_F^*Ti_F(y)$.  

Suppose, now, that both $r_0$ and $F_0$ are
sufficiently large so that there exists $y_0 \in F_0$ with $\|y_0\|
\leq r_0$ and $y^*_0 \in i_{F_0}^*Ti_{F_0}(y_0).$   Fix $F \supset
F_0$; for each $r \geq r_0$ we will show that the following set $H_r$
is bounded:  Define 
$H_r$ to be the set of all $(x_{_{F,r}}, x_{_{F,r}}^*) \in F \times
F^*$ such that  $$x_{_{F,r}}^* \in i_F^*J(x_{_{F,r}})\ {\rm and}\ 
\langle -x_{_{F,r}}^* - y^*, x_{_{F,r}} - y \rangle \geq 0\ {\rm 
whenever}\ \|y \| \leq r\ {\rm and}\ y^* \in i_F^*Ti_F(y).$$
  For any such $(x_{_{F,r}},
x_{_{F,r}}^*)$ we have  $${1\over 2}\|x_{_{F,r}}\|^2 + {1\over
2}\|x_{_{F,r}}^*\|^2 = \langle x_{_{F,r}}^*,x_{_{F,r}}\rangle \leq
\langle x_{_{F,r}}^*, y_0\rangle - \langle y_0^*, x_{_{F,r}}\rangle +
\langle y_0^*, y_0\rangle \leq$$ $$\leq \|x_{_{F,r}}^*\|\cdot \|y_0\| +
\|y_0^*\|\cdot \|x_{_{F,r}}\| +  \langle y_0^*, y_0\rangle.$$ The
subset of the plane where a positive quadratic function is dominated
by a linear function is necessarily bounded, so there exists an upper
bound on each of the sets $\{ \|x_{_{F,r}}\|\}$ and
$\{\|x_{_{F,r}}^*\|\}$.  That is, each of the sets $H_r$ is bounded. 
Moreover, each of them is closed (in the product of the norm
topologies), since $i_F^*Ji_F$ is easily seen to have closed graph and
the function $(x_{_{F,r}}, x_{_{F,r}}^*) \to  \langle -x_{_{F,r}}^* -
y^*, x_{_{F,r}} - y \rangle$ is continuous, for each $(y, y^*)  \in F
\times F^*$.  Clearly, $H_r \supset H_{r'}$ whenever $r' > r > 0$, so
for increasing $r$, the $H_r$'s form a decreasing family of
nonempty compact sets and they therefore have nonempty intersection. 
This shows that there exists $x_{_F} \in F$ and $x_{_F}^* \in
i_F^*Jx_{_F}$ such that   $$\langle -x_{_F}^* - y^*, x_{_F} - y
\rangle \geq 0\ \hbox {\rm whenever}\ y \in F,\  \hbox {\rm and}\ y^*
\in i_F^*Ti_F(y).\eqno (1)$$  
Note that any Hahn-Banach
extension of $x_{_{F}}^*$ from $F$ to all of $E$ is in 
$J(x_{_{F}})$, so we can assume that $x_{_{F}}^* \in
J(x_{_{F}}) \subset E^*$.  Since the nets $(x_{_F})$ and $(x_{_F}^*)$
are bounded,  and since we can regard the $x_{_F}$'s as elements of
$E^{**}$, we see that there exists a subnet (call it $(x_{_F},
x_{_F}^*)$) in $E^{**} \times E^*$ converging to an element $(x^{**},
x^*) \in E^{**} \times E^*$ in the $\sigma (E^{**}, E^*) \times \sigma
(E^*, E)$ topology.  We want to show that $(x^{**}, -x^*)$ is
monotonically related to $G(T) \subset E^{**} \times E^*$, that is,
$$\langle x^{**} - \hat y, -x^* - y^*\rangle \geq 0 \ \hbox {\rm
whenever}\ (y,y^*) \in G(T). \eqno (2)$$ To see this, note that (using
the weak* lower semicontinuity of both $j^*$ and $j^{**}$) $$\langle
x^{**}, x^*\rangle \leq j^{**}(x^{**}) + j^*(x^*) \leq \liminf
[j(x_{_F}) + j^*(x_{_F}^*)] = \liminf \langle x_{_F}^*,
x_{_F}\rangle,\eqno (3)$$ while (1) implies that, for all $(y, y^*)
\in G(T)$,
  $$ \limsup
\langle x_{_F}^*, x_{_F}\rangle \leq \langle y^*, y\rangle + \langle
x^*, y\rangle - \langle x^{**}, y^*\rangle; \eqno (4)$$ together, these
yield (2).  Since $T$ is assumed to be of type (D), there exists a
net $(y_\alpha, -y_\alpha^*)$ in $G(T)$ such that $(y_\alpha)$ is
bounded, converges to $x^{**}$ in the $\sigma (E^{**}, E^*)$ topology
and $\|y_\alpha^* - x^*\| \to 0$.  This fact, applied to (4), shows
that  $$\limsup \langle x_{_F}^*, x_{_F}\rangle \leq \langle x^{**},
x^*\rangle + \langle x^{**}, x^*\rangle - \langle x^{**}, x^*\rangle =
\langle x^{**}, x^*\rangle.$$ Now, from (3), $$\langle x^{**},
x^*\rangle \leq j^{**}(x^{**}) + j^*(x^*) \leq
 \liminf \langle x_{_F}^*, x_{_F}\rangle \leq \limsup \langle
x_{_F}^*, x_{_F}\rangle \leq \langle x^{**}, x^*\rangle,$$
which shows that $x^{**} \in \partial j^*(x^*)$. Thus, $-x^* \in
\overline T(x^{**}) $ and $x^* \in \partial (j^*)^{-1}(x^{**})$, which
completes the proof.

\proclaim Corollary 3.7.  If $E$ is reflexive and $T\colon E \to 2^{E^*}$ is
maximal monotone, then $R(T + \lambda J) = E^*$ for every $\lambda > 0$.

\noindent {\bf Proof.}  It follows directly from the definitions that for
a reflexive space $E$, one has $\partial (j^*)^{-1} = J$ and, as has been noted
earlier, $\overline T = T$, so the corollary is immediate.

\proclaim Theorem 3.8.  If $T$ is of maximal monotone of type (D), then
$\overline {R(T)}$ is convex.

\noindent Before proving this, we need the fact that if $T$ is maximal monotone
of type (D) and $x^* \in {\rm co}R(T)$, then there exists $x \in E$ such that
$$\sup_{(y^{**}, y^*) \in G(\overline T)} \langle y^{**} - \hat x,
x^* - y^*\rangle < \infty.$$

\noindent This follows in a straightforward way from the following lemma and the
definition of type (D).  Note that if $T$ is of type
(D), then $R(\overline T) \subset \overline {R(T)}.$

\proclaim Lemma 3.9. Suppose that $E$ and $F$ are linear
spaces in duality and that $T\colon E \to 2^F$ is monotone.  If $x^* \in {\rm
co} R(T)$, then there exists $x \in {\rm co} D(T)$ such that $$\sup_{(y, y^*) \in
G(T)}\langle y^* - x^*, x - y\rangle < \infty. $$

\noindent {\bf Proof.}  Suppose that $x^* = \sum t_i x_i^*$ where $t_i \geq
0$, $\sum t_i = 1$ and $x_i^* \in R(T)$, so there exist $x_i \in E$ such
that $x_i^* \in T(x_i)$.  Take $x = \sum t_ix_i$; then for any $(y, y^*)
\in G(T)$,
 
$$\langle y^* - x^*, x - y\rangle = \langle y^* - \sum t_ix_i^*, \sum
t_jx_j - y\rangle = \sum_{i, j}t_it_j\langle y^* - x_i^*, x_j - y\rangle =$$

$$= \sum_{i, j}t_it_j\langle y^* - x_j^*, x_j - y\rangle + \sum_{i,
j}t_it_j\langle x_j^* - x_i^*, x_j - y\rangle \leq \sum_{i, j}t_it_j\langle x_j^* -
x_i^*, x_j - y\rangle =$$

$$= \sum_{i< j}t_it_j\langle x_j^* - x_i^*, x_j - x_i\rangle ,$$ which proves
the lemma, since the last term does not depend on $y$ or $y^*$.

 \medskip 
\noindent {\bf Proof of Theorem 3.8.}  Note that it suffices to show that
${\rm co}R(T) \subset \overline {R(T)}$, since this implies that
$\overline {{\rm co}}R(T) \subset \overline { R(T)} \subset \overline
{{\rm co}}R(T)$.  Suppose, then, that $x^* \in {\rm
co}R(T)$.  By Theorem 3.6, for each $\lambda > 0$ there exists
$y_\lambda ^* \in E^*$, $x_\lambda ^{**} \in j^*(y_\lambda^*)$ and
$z_\lambda ^* \in \overline T(x_\lambda ^{**})$ such that $x^* =
\lambda y_\lambda ^* + z_\lambda ^*$.  By the foregoing remark, there
exists $x \in E$ such that $$\langle x_\lambda ^{**} - \hat x, x^* -
z_\lambda ^*\rangle \equiv \lambda \langle x_\lambda ^{**} - \hat x,
y_\lambda ^*\rangle
 $$
is bounded above for all $\lambda > 0$.  It follows that for some $M
> 0$ (and all $\lambda > 0$),
 $$\lambda \| y_\alpha ^*\|^2 \leq \lambda \| 
x_\lambda^{**} \|^2 + \lambda \| y_\lambda ^*\|^2 = 2\lambda
\langle x_\lambda ^{**}, y_\lambda ^*\rangle \leq M + 2\lambda \langle
y_\lambda ^*, x\rangle \leq M + 2\lambda \|y_\lambda^*\| \cdot \|x\|.$$
>From this we see that $\lambda y_\lambda ^* \to 0$ as $\lambda \to
0$; indeed, if there were a sequence $\lambda_n \to 0$ such that
$\| \lambda_n y_{\lambda_n}^*\|$ were bounded away from $0$, then we
would necessarily have $\|y_{\lambda_n}^*\| \to \infty$ and dividing
both sides of the inequality above by $\lambda_n \|y_{\lambda_n}^*\|$
would lead to a contradiction.  Thus, $x^* - z_\lambda ^* = \lambda
y_\lambda ^* \to 0 $; since $z_\lambda ^* \in R(\overline T)$, this
shows that $x^* \in \overline {R(\overline T)} \subset \overline
{R(T)}$.
\medskip \noindent (Since $R(T) \subset R(\overline T)$, their closures
are in fact equal.)

\proclaim Corollary 3.10.  If $E$ is reflexive and $T\colon E \to 2^{E^*}$ is
maximal monotone, then both $\overline {R(T)}$ and $\overline {D(T)}$ are convex.

\noindent {\bf Proof.}  By Example 3.2(b), if $T$ is maximal monotone, then it
is maximal monotone of type (D), so it follows from Theorem 3.8 that $\overline
{R(T)}$ is convex.  By applying this result to the maximal monotone operator
$T^{-1}$ we obtain convexity of $\overline {R(T^{-1})} \equiv \overline {D(T)}$.

\proclaim Definition 3.11.  An operator $T\colon E \to 2^{E^*}$ is said to be
{\it coercive} provided $D(T)$ is bounded or there exists a function $c\colon
{\bf R}^+ \to {\bf R}$ such that $c(r) \to \infty$ when $r \to \infty$ and
$\langle x^*, x\rangle \geq c(\|x\|)\cdot \|x\|$ for each $(x, x^*) \in G(T)$.

\medskip \noindent {\bf Remark.}  It is easily 
verified that if $D(T)$ is unbounded, then $T$ is coercive if and only if for
every $M > 0$ there exists $r > 0$ such that 
$${\langle x^*, x\rangle \over \|x\|} \geq M\quad {\rm
whenever}\quad \|x\| \geq r\ {\rm and}\ x^* \in T(x).$$
Indeed, if this holds, take $c(r) = \inf \{ {\langle x^*, x\rangle \over
\|x\|}\colon\quad \|x\| \geq r \ {\rm and} \ x^* \in T(x)\}$.

\medskip \noindent {\bf Examples 3.12.}  (a)  The
duality mapping $J$ is an obvious example of a coercive
operator, since
$\langle x^*, x\rangle = \|x\|^2$ whenever $x^* \in J(x)$.

\smallskip (b)  If $T$ is a positive linear operator and $\lambda > 0$, then $T
+ \lambda J$ is coercive:  If $x \in E$ and $x^* \in (T + \lambda J)(x)$, then
$x^* = T(x) + \lambda z^*$ for some $z^* \in J(x)$ and hence $\langle x^*,
x\rangle = \langle T(x), x\rangle + \lambda \langle z^*, x\rangle \geq \lambda
\|x\|^2$.  From the Remark following Problem 2.20, it follows that $T
+ \lambda J$ is also maximal monotone.

\smallskip (c)  Recall that a function $\varphi\colon {\bf R}\to {\bf R}$ is a
monotone operator if and only if it is nondecreasing.  It is easily seen that
$\varphi$ is coercive if and only if $\varphi (t) \to \pm \infty$ as $t \to \pm
\infty$.

\proclaim Exercise 3.13.  Show that if $T$ is coercive, then so is $\overline
T\colon E^{**} \to 2^{E^*}$.

\proclaim Theorem 3.14 (Gossez).  Suppose that $T$ is a coercive maximal
monotone operator of type (D).  Then $R(\overline T) = E^*$ and hence $\overline
{R(T)} = E^*$.

\noindent {\bf Proof.}  It is clear from the definition of type (D) that one
always has $R(\overline T) = \overline {R(T)}$.  Suppose, then, that $x^* \in
E^*$.  By Theorem 3.6, for each $\lambda > 0$ there exist $y_\lambda^* \in
E^*$, $x_\lambda^{**} \in \partial j^*(y_\lambda^*)$ and $z_\lambda^* \in
\overline T(x_\lambda^{**})$ such that $x^* = \lambda y_\lambda^* +
z_\lambda^*$.  We have
$$\langle x_\lambda^{**}, x^*\rangle = \lambda \langle x_\lambda^{**},
y_\lambda^*\rangle + \langle x_\lambda^{**}, z_\lambda^*\rangle =
\lambda \|x_\lambda^{**}\|^2 + \langle x_\lambda^{**}, z_\lambda^*\rangle.$$
Since (Exercise 3.13) $\overline T$ is coercive, if $\{\|
x_\lambda^{**}\| \}$ were unbounded as $\lambda \to 0$, the right side of 
$$\displaystyle \|x^*\| \geq {\langle
x_\lambda^{**}, x^*\rangle \over \|x_\lambda^{**}\|} = 
\lambda \|x_\lambda^{**}\| + {\langle x_\lambda^{**}, z_\lambda^*\rangle \over 
\|x_\lambda^{**}\|}$$
would be unbounded, an impossibility.  Thus, the bounded net $\{x_\lambda^{**}\}$
has a subnet (call it $\{x_\lambda^{**}\}$) converging in the $\sigma (E^{**},
E^*) $ topology to an element $x^{**} \in E^{**}$.  We will show that $x^* \in
\overline T(x^{**})$ by showing that $(x^{**}, x^*)$ is monotonically related to
$G(T)$.  Suppose, then, that $(u, u^*) \in G(T)$.  Since $x_\lambda^* = x^* -
\lambda y_\lambda^*$, we have 
$$0 \leq \langle x_\lambda^{**} - \hat u, x_\lambda^* - u^*\rangle
= \langle x_\lambda^{**} - \hat u, x^* - u^*\rangle - \lambda \langle
x_\lambda^{**} - \hat u, y_\lambda^*\rangle.$$
Recall that $x_\lambda^{**} \in \partial j^*(y_\lambda^*)$ implies boundedness of
$\|y_\lambda ^*\| = \|x_\lambda^{**}\|$, so the second term on the right
converges to $0$ as $\lambda \to 0$, yielding $0 \leq \langle x^{**} -
\hat u, x^* - u^*\rangle$.

\proclaim Corollary 3.15.  If $E$ is reflexive and $T$ is a coercive maximal
monotone operator on $E$, then $R(T) = E^*$.

\noindent {\bf Proof.}  Simply use the fact that reflexivity implies that $T =
\overline T$.

\medskip \noindent The fact (Proposition 2.17(c)) that $E$ is reflexive if $R(J)
= E^*$ shows that one cannot omit reflexivity from this result.  

\bigskip  \noindent In order that $\overline {R(T)}$ be convex it is {\it
not} necessary for $T$ to be maximal monotone of type (D).  This was shown by
Gossez [Go$_{2,3}$] with the help of the following example.

\bigskip \noindent {\bf Example 3.16.}  Let $A$ be defined on the nonreflexive
space $\ell_1$ as follows:  For each $x = (x_k) \in \ell_1$, let $\{A(x)_n\}$
be the $\ell_\infty$ sequence defined by
 $$(Ax)_n = -\sum_{k  < n}x_k + \sum_{k > n}x_k.$$ It is not hard to verify
that $A$ is bounded, linear and antisymmetric (that is, $\langle Ax, y\rangle
= -\langle Ay, x\rangle$ for all $x, y \in \ell_1$) hence is monotone,
satisfying $\langle Ax, x\rangle = 0$ for all $x$.  This latter means that,
in particular, $A$ is a positive operator, hence it is maximal monotone.  The
range $R(A)$ of $A$ is a linear subspace (hence is convex and has convex
closure) which is properly contained in the proper closed subspace $c$ of
$\ell_\infty$ consisting of all convergent sequences.  (Indeed, $\lim_{n \to
\infty}(Ax)_n = -\sum_{k = 1}^\infty x_k$.)  

\medskip \noindent Gossez [Go$_2$] uses the operator $A$ by showing that there
exists $\lambda > 0$ such that $R(A + \lambda J)$ is not dense in
$\ell_\infty$.  This shows that $A + \lambda J$ is not of type (D), in view of
Theorem 3.14 and the fact that $A + \lambda J$ is maximal monotone and coercive
(Example 3.12(b)).  Thus {\it not all maximal monotone operators (not even the
coercive ones) are of type (D).}  Subsequently, he showed [Go$_3$] that the fact
that $R(A + \lambda J)$ is not dense in $\ell_\infty$ implies that its closure
is not convex, that is, {\it there exists a coercive maximal monotone
operator $T$ on $\ell_1$ such that $\overline {R(T)}$ is not convex.}

\bigskip \bigskip

\magnification = \magstep1
\noindent {\bf 4.  Locally  maximal monotone operators.}

\bigskip  As we have seen, some of the nice properties of maximal monotone
operators on reflexive spaces fail to hold in general, but {\it are} valid
for the subclass of maximal monotone operators of type (D).  In this section
we introduce another subclass which shares some of the same properties.

 \proclaim  Definition 4.1.  A
set-valued mapping $T\colon E \to 2^{E^*}$ is said to be {\it locally maximal
monotone} if, for each norm-open convex subset $U \subset E^*$ which
intersects $R(T)$, the restriction of the inverse operator $T^{-1}$ to $U$ is
maximal monotone in $U$.  The latter means that the graph $G((T^{-1})|_U)
\subset U \times E$ is a maximal monotone subset of $U\times E$.  

\noindent The
``working definition'' of this property is the following:  If $U$ is an open
convex subset of $E^*$ which intersects $R(T)$ and if $(x, x^*) \in E\times U$
is montonically related to each $(y, y^*) \in G(T) \cap (E \times U)$, then 
$(x, x^*) \in G(T)$.
\medskip
\noindent  It is clear (take $U = E^*$) that every locally maximal monotone
operator is maximal monotone.
\bigskip

The locally maximal monotone operators were introduced in [F--P$_1$] because
they are the precise class for which a certain approximation scheme is
valid.  While their exact position within the class of all maximal monotone
operators is still unclear, some important properties are known.

 \proclaim
Proposition 4.2.   \smallskip (i)  If $T$ is locally maximal monotone, then
$\overline {R(T)}$ is convex.
\smallskip (ii)  If $f$ is a proper lower semicontinuous convex function on
$E$, then $\partial f$ is locally maximal monotone.

\noindent  The proof for (i) may be found in [F--P$_1$].  Property
(ii), which is a nontrivial extension of Rockafellar's maximality theorem
(Theorem 2.15), was proved by S. Simons [Si$_2$]; see, also, [Si$_3$].  In
order to see that maximal monotone operators in reflexive spaces are locally
maximal monotone, we first reformulate the definition.

\proclaim Proposition 4.3. A monotone 
operator $T$ on $E$ is locally maximal
monotone if and only if it satisfies the following condition: For any 
weak* closed convex and bounded subset $C$ of $E^*$ such
that $ R(T) \cap {\rm int}\thinspace C \not= \emptyset$ and for each $x \in 
E$ 
and
$x^* \in  {\rm int}\thinspace C$ with $x^* \notin T(x)$, there exists $z \in E 
$ 
and 
$z^* \in T(z)\cap C$ such that $\langle x^* - z^*, x - z\rangle < 0$.

\noindent
\bf Proof. \rm In one direction, if $T$ is locally maximal monotone 
and $C$ is given, let $U = {\rm int}\thinspace C$.  In the other direction,
 if $U$ is
open and convex in $E^*$, if $u \in E$ and $x \in E$ with $~{u^* \in 
T(u)\cap U}$ and $x^* \in U$
but $x^* \notin T(x)$, then there exists $\epsilon > 0$ such that $u^* +
\epsilon B^* \subset U$ and $x^* + \epsilon B^* \subset U$.  By 
convexity,
$C \equiv [u^*, x^*] + \epsilon B^* $ is a weak* closed, convex
and bounded subset of $U$ which can be used to verify that $U$ has 
the required property.

\bigskip The only use of reflexivity in the next proposition is an
application of Theorem 2.19 (on the sum of two maximal monotone operators in
a reflexive space).

\bigskip \noindent
\proclaim  Proposition 4.4. If $E$ is reflexive and $T$ is maximal
monotone on $E$, then it is locally maximal monotone.

\noindent \bf Proof. \rm Suppose that $C$ is weak* closed and convex and that
${\rm int}\thinspace C \cap R(T) \neq \emptyset$.  Suppose also that $x^* \in {\rm
int}\thinspace C$ but $x^* \notin T(x)$.  Let $T_1$ denote the inverse $T^{-1}$
of $T$ and let $T_2 = \partial \delta_C$.  Since ${\rm int}\thinspace D(T_2) =
{\rm int}\thinspace C$, these are maximal monotone operators from $E^*$ into $E$
for which $D(T_1) \cap {\rm int}\thinspace D(T_2) \neq \emptyset$; by Theorem
2.19, their sum $T_1 + T_2$ is maximal monotone.
Now
$x^* \notin T(x)$ implies that $x \notin T_1(x^*)$, 
and since $T_2(x^*) = \{0\}$,
we see that $x \notin T_1(x^*) + T_2(x^*)$.  By maximality of 
$T_1 + T_2$, there exists $z^* \in D(T_1) \cap D(T_2) \equiv
R(T) \cap C$ and $z \in (T_1 + T_2)(z^*)$ such that 
$\langle x^* - z^*, x - z\rangle < 0$.  We can write $z = u + v$, where $u \in
T_1(z^*)$ (that is, $z^* \in T(u)$) and $v \in T_2(z^*)$.  The latter means that
$\langle z^* - w^*, v\rangle \geq 0$ for all $w^* \in C$.  We have thus produced
$z^* \in T(u) \cap C$ such that 

$$0 > \langle x^* - z^*, x - u\rangle - \langle x^* - z^*, v\rangle \geq \langle
x^* - z^*, x - u\rangle,$$
showing that $T$ satisfies the condition in Proposition 4.3  and is therefore
locally  maximal monotone.

\medskip  Recall Gossez's Example 3.15 of the linear maximal monotone operator
$A\colon \ell_1 \to \ell_\infty$.  Its interest in this context is the
fact that, even though $\overline {R(A)}$ is linear (hence convex), $A$ is not
locally maximal monotone, so {\it not every maximal monotone operator is
locally maximal monotone.}

\proclaim Example 4.5.  The operator $A$ is not locally maximal monotone.

\noindent {\bf Proof.} Let $e = (1, 0, 0, \dots)$, considered as an element of
either $\ell_1$ or $\ell_\infty$, and let $$z = (-{1\over 2}, {1\over 2^3},
{1\over 2^4}, {1\over 2^5}, \dots) \in \ell_1.$$  Some computations using the
definition show that $(Az)_1 = {1\over 4}$ while for $n \geq 2$, $(Az)_n =
{1\over 4} + {1\over 2^n} + {1\over 2^{n + 1}}$.  Moreover, $e - Ae = (1, 1,
1,\dots)$ and $\|e - Az\|_\infty = {3\over 4}$, so if $U$ is the open unit
ball in $\ell_\infty$, then $x^* \equiv e - Az \in U$.  Let $x = e - z$.  If
$u \in \ell_1$ and $Au \in U$, then $\lim_{n \to \infty} |(Au)_n| = |\sum
_{k = 1}^\infty u_k| \leq 1$ and hence $$\langle x^* - Ax, u\rangle = \langle
e - Ae, u\rangle = \sum _{k = 1}^\infty u_k \leq 1,$$ while $$\langle x^*,
x\rangle = \langle e, e\rangle - \langle Az, e\rangle - \langle e, z\rangle =
1 - (Az)_1 - z_1 = 1 - {1\over 4} + {1\over 2} > 1.$$ Thus, $x^* \neq Ax $
even though $$\langle x^* - Au, x - u\rangle = \langle x^*, x\rangle -
\langle Au, x\rangle - \langle x^*, u\rangle = \langle x^*, x\rangle +
\langle Ax, u \rangle - \langle x^*, u\rangle \geq 0$$ whenever $Au \in U$,
contradicting the definition of locally maximal monotone. 

\bigskip {\it We still do not know whether the class of maximal monotone
operators of type (D) is actually different from the class of locally maximal
monotone operators.}  To see that for {\it coercive} operators, the
former class is contained in the latter, we need two preliminary results.  The
proof of the following identity consists of an elementary but tedious computation.

 \proclaim Proposition 4.6.  If  $u,v,x\in E$,
$u^*,v^*,x^*\in E^*$ and $\lambda \in [0, 1]$, then   $$\eqalign{\langle\lambda
u^*+(1-\lambda)v^*-x^*,\lambda
u+(1-\lambda)v-x\rangle  = \qquad \qquad \qquad \cr = \lambda\langle
u^*-x^*,u-x\rangle +(1-\lambda)\langle
v^*-x^*,v-x\rangle- \lambda(1-\lambda)\langle
u^*-v^*,u-v\rangle \cr} \eqno (4.1)$$

\proclaim Lemma 4.7.   Suppose that $T$ is a maximal monotone
operator, that $U$ is an open subset of $E^*$
and that $z^*
\in U\setminus T(z)$ is such that $\langle x^*-z^*,x-z\rangle\ge 0$ for all
$x^*\in T(x)\cap U$. Then there
exist $b\in E$, $b^*\in U$ and $r>0$ such that 
for all $x^*\in T(x)\cap U$,
$$\langle x^*-b^*,x-b\rangle\ge r.$$

\noindent {\bf Proof. } Since $T$ is maximal monotone there exists
$y^*\in T(y)$ such that $\langle y^*-z^*,y-z\rangle<
0$. Let $1>\lambda>0$ be such that $b^*:=\lambda z^*+
(1-\lambda)y^*\in U$ and let $b =\lambda z+
(1-\lambda)y$. Then, using the identity (4.1), for all $x^*\in T(x)\cap U$ we
have $$\eqalign{\langle x^*-b^*,x-b\rangle&=
\lambda\langle x^*-z^*,x-z\rangle+
(1-\lambda)\langle x^*-y^*,x-y\rangle -
\lambda(1-\lambda)\langle z^*-y^*,z-y\rangle\cr
&\ge -\lambda(1-\lambda)\langle z^*-y^*,z-y\rangle>0
}$$
so we may set $r =-\lambda(1-\lambda)\langle
z^*-y^*,z-y\rangle>0$.

\proclaim Theorem 4.8.  Suppose that $T$ is a maximal monotone operator such
that either (i) $R(T) = E^*$ or (ii) $\overline { R(T) } =
E^*$ and $T$ is coercive.  Then $T$ is locally maximal monotone.

\noindent {\bf Proof.}  Suppose, first, that $R(T) = E^*$, that $U \subset
E^*$ is open and convex and that $z \in E$, $z^* \in U$ are such that
$\langle z^* - x^*, z - x\rangle \geq 0$ for all $x \in E$ such that $x^* \in
T(x) \cap U$.  If $z^* \notin T(z)$, then there would exist $b \in E$, $b^* \in
U$ and $r > 0$ as in Lemma 4.7.  Since $b^* \in R(T)$ by hypothesis, there
exists $x \in E$ such that $b^* \in T(x) \cap U$ and hence by Lemma 4.7,
$\langle b^* - b^*, x - b\rangle \geq r > 0$, a contradiction.

\smallskip \noindent Suppose, next, that $R(T)$ is dense in $E^*$ and that $T$
is coercive.  If $T$ were not locally maximal monotone, we could find an open
convex subset $U \subset E^*$ with $U \cap R(T) \neq \emptyset$ and elements
$z \in E$ and $z^* \in U\backslash T(z)$ such that $\langle z^* - x^*, z -
x\rangle \geq 0$ whenever $x \in E$ and $x^* \in T(x) \cap U$.  Choose $b$,
$b^*$ and $r$ as in Lemma 4.7.  Since $R(T)$ is dense in $E^*$, we can find
$x_n \in E$ and $x_n^* \in T(x_n)$ such that $\|b^* - x_n^*\| \to 0$.  But then
for all sufficiently large $n$, we would have $x_n^* \in U$ and hence
$$r \leq \langle x_n^* - b^*, x_n - b\rangle \leq \|x_n^* - b^*\|  \|x_n -
b\|,$$
which would imply that $\|x_n - b\| \to \infty$.  Coercivity would then imply
that $\|x_n^*\| \to \infty$, again a contradiction.

\proclaim Corollary 4.9.  If $T$ is maximal monotone, coercive and of 
type (D), then it is locally maximal monotone.

\noindent {\bf Proof.}  Recall that by Theorem 3.14, the fact that $T$
is coercive and maximal monotone of type (D) implies that $\overline {R(T)} =
E^*$.  

\bigskip Recall Problem 2.20:  If $E$ is a nonreflexive Banach space and $S$
and $T$ are maximal monotone operators such that $D(T) \cap {\rm int}\thinspace
D(S) \neq \emptyset$; is $S + T$ necessarily maximal?  What about the special
case when $S$ is the subdifferential of the indicator function $\delta_C$ of a
closed convex set $C$ for which ${\rm int}\thinspace C \cap D(T) \neq \emptyset$?

\smallskip \noindent  It is not unreasonable to ask whether these questions have
affirmative answers when the maximal monotone operators are of type (D), or are
locally maximal monotone. 

\bigskip \bigskip

\baselineskip=15truept
\vsize=8.5truein 
\hsize=6.5truein 
\parskip=18truept

\def\inch{{{\vrule width.73in depth0.15pt height0.15pt}}}

\centerline{\bf References}

\item{[B-F]}	J. Borwein and S. P. Fitzpatrick, ``Local
boundedness of monotone  operators under minimal hypotheses'',
{\sl Bull. Australian Math. Soc.} {\bf 39} (1989), 439-441.

\item{[Bre]}H.~Brezis, ``Op\'erateurs Maximaux Monotones
et semi-groupes de contractions dans les espaces de
Hilbert'', Math. Studies 5, North-Holland American
Elsevier (1973).

\item {[B-C-P]} H. Brezis, M. Crandall and A. Pazy, ``Perturbations of nonlinear
maximal monotone set'', {\sl Comm. Pure. Appl. Math.} {\bf 23} (1970), 123-144.

\item{[Bro]}F. E. Browder, ``The fixed point theory of
multi--valued mappings in topological vector spaces'', {\sl Math.
Annalen} {\bf 177} (1968), 283--301.

\item{[De]}K. Deimling, ``Nonlinear Functional
Analysis'', Springer-Verlag (1985).

\item{[D-G-Z]} R. Deville, G. Godefroy and V. Zizler,
``Renormings and Smoothness in Banach Spaces'', Monographs and
Surveys in Pure and Appl. Math., Longman (1993).

\item{[Di]}J. Diestel, ``Geometry of Banach
space---selected topics'', Lecture Notes in Math., Nr. 485,
Springer-Verlag (1975).

\item{[F-$\hbox{\rm P}_1$]}S. P. Fitzpatrick and R. R.
Phelps, ``Bounded approximants to monotone operators on
Banach spaces'', {\sl Ann. Inst. Henri Poincar\'e, Analyse
non lin\'eaire} {\bf 9} (1992), 573--595.

\item{[F-$\hbox{\rm P}_2$]}\inch, ``Maximal monotone
operators on nonreflexive Banach spaces'' (in preparation).

\item{[Fr]}J. Franklin, ``Methods of Mathematical Economics'',
Springer--Verlag, Undergraduate Texts in Math. (1980).

\item{[$\hbox{\rm Go}_1$]}J.-P. Gossez, ``Op\'erateurs
monotones non lin\'eaires dans les espaces de Banach non
r\'eflexifs'', {\sl J. Math. Anal. Appl.} {\bf 34} (1971),
371--395.

\item{[$\hbox{\rm Go}_2$]}\inch, ``On the range of a
coercive maximal monotone operator in a nonreflexive
Banach space'', {\sl Proc. Amer. Math. Soc.} {\bf 35}
(1972), 88--92.

\item{[$\hbox{\rm Go}_3$]}\inch, ``On a convexity property
of the range of a maximal monotone operator'', {\sl Proc.
Amer. Math. Soc.} {\bf 55} (1976), 359--360.

\item{[$\hbox{\rm Go}_4$]}\inch, ``On the extensions to
the bidual of a maximal monotone operator'', {\sl Proc.
Amer. Math. Soc.} {\bf 62} (1977), 67--71.

\item{[Go$_5$]}\inch, private communication.

\item{[Pa-Sb]}D. Pascali and S. Sburlan, ``Nonlinear
mappings of monotone type'', Ed. Acad., Bucarest, Rom.,
Sijthoff \& Noordhoff Internat. Publ., Alphen aan den
Rijn, Netherlands (1978).

\item{[Ph]}R. R. Phelps, ``Convex Functions, Monotone
Operators and Differentiability'',  Lecture Notes
Math. Nr.1364, Springer-Verlag (1989).  Second Edition
(1993).

 \item{[Re]}S. Reich,
``The range of sums of accretive and monotone operators'', {\sl J. Math. Anal.
Appl.} {\bf 68} (1979), 310--317.

\item{[R-V]}A. W. Roberts and D. E. Varberg, ``Convex
Functions'', Academic Press, New York-San Francisco-London
(1973).

\item{[$\hbox{\rm Ro}_1$]}R. T. Rockafellar, ``On the
maximality of sums of nonlinear monotone operators'', {\sl
Trans. Amer. Math. Soc.} {\bf 149} (1970), 75--88.

\item{[$\hbox{\rm Ro}_2$]}\inch, ``Local boundedness of
nonlinear, monotone operators'', {\sl Mich. Math. J.} {\bf
16} (1969), 397--407.

\item{[$\hbox{\rm Ro}_3$]}\inch, ``On the maximal
monotonicity of subdifferential mappings'', {\sl Pacific
J. Math.} {\bf 44} (1970), 209--216.

\item{[$\hbox{\rm Ro}_4$]}\inch, ``On the virtual convexity of
the domain and range of a nonlinear maximal monotone
operator'', {\sl Math. Ann.} {\bf 185} (1970), 81--90.

\item{[$\hbox{\rm Si}_1$]}S. Simons, ``The least slope of
a convex function and the maximal monotonicity of its
subdifferential'', {\sl J. Optimization Theory and
Applications} {\bf 71} (1991), 127--136.

\item{[$\hbox{\rm Si}_2$]}\inch, ``Subdifferentials are
locally maximal monotone'', {\sl Bull. Australian Math.
Soc.} {\bf 47} (1993), 465--471.

\item{[$\hbox{\rm Si}_3$]}\inch, ``Les d\'eriv\'ees
directionelles et la monotonicit\'e des
sous-diff\'erentiels'', S\'em. d'Initiation \`a l'Analyse
(S\'em. Choquet), Paris (to appear).

 \item{[Ze]}E. Zeidler, ``Nonlinear
Functional Analysis and its Applications'', Vol. II/A Linear
Monotone Operators, Vol. II/B Nonlinear Monotone Operators,
Springer-Verlag (1985).

\end